\documentclass{article}
\usepackage{cancel}
\usepackage{graphicx}
\usepackage{subcaption}
\usepackage{geometry}
\usepackage{adjustbox}
\usepackage{bm}
\usepackage{algorithmicx}
\usepackage[ruled,vlined]{algorithm2e}
\usepackage{algpseudocode,float}
\algblock{ParFor}{EndParFor}

\usepackage{mathrsfs}  
\usepackage{sidecap}
\usepackage{amsfonts}
\usepackage{amssymb}
\usepackage{amsmath}
\usepackage{amsthm}
\usepackage{comment}
\usepackage{xcolor}
\usepackage{relsize}
\usepackage{multirow}
\usepackage{multicol}
\usepackage{diagbox}
\usepackage{overpic}
\usepackage{authblk}

\newcommand{\mC}{{\mathsf C}}

\newcommand{\mT}{{\mathsf T}}
\newcommand{\mL}{{\mathsf L}}
\newcommand{\mS}{{\mathsf S}}
\newcommand{\mA}{{\mathsf A}}
\newcommand{\mB}{{\mathsf B}}
\newcommand{\mD}{{\mathsf D}}

\newcommand{\mI}{{\mathsf I}}

\newcommand{\mW}{{\mathsf W}}
\newcommand{\mSigma}{{\mathsf \Sigma}}

\newcommand{\R}{{\mathbb R}}
\newtheorem{theorem}{Theorem}

\usepackage{cleveref}
\usepackage{autonum}
\date{}
\begin{document}
\title{Overcomplete representation in a hierarchical Bayesian framework}
%
%
\author[1]{Monica Pragliola\thanks{monica.pragliola2@unibo.it}}
\author[2]{Daniela Calvetti\thanks{dxc57@case.edu}}
\author[2]{Erkki Somersalo\thanks{ejs49@case.edu}}
\affil[1]{Department of Mathematics, University of Bologna, Italy}
\affil[2]{Case Western Reserve University, Department of Mathematics, Applied Mathematics and Statistics, 10900 Euclid Avenue, Cleveland, OH 4410, USA}

%
%
%
\maketitle              

\begin{abstract}
	A common task in inverse problems and imaging is finding a solution that is sparse, in the sense that most of its components vanish. In the framework of compressed sensing, general results guaranteeing exact recovery have been proven. In practice, sparse solutions are often computed combining $\ell_1$-penalized least squares optimization with an appropriate numerical scheme to accomplish the task. A computationally efficient alternative for finding sparse solutions to linear inverse problems is provided by Bayesian hierarchical models, in which the sparsity is encoded by defining a conditionally Gaussian prior model with the prior parameter obeying a generalized gamma distribution. An iterative alternating sequential (IAS) algorithm has been demonstrated to lead to a computationally efficient scheme, and combined with Krylov subspace iterations with an early termination condition, the approach is particularly well suited for large scale problems. Here the Bayesian approach to sparsity is extended to problems whose solution allows a sparse coding in an overcomplete system such as composite frames. It is shown that among the multiple possible representations of the unknown, the IAS algorithm, and in particular, a hybrid version of it, is effectively identifying the most sparse solution. Computed examples show that the method is particularly well suited not only for traditional imaging applications but also for dictionary learning problems in the framework of machine learning.
\end{abstract}

\section{Introduction}
\label{sec:intro}

Sparsity promoting methods and algorithms for inverse problems and imaging applications have been extensively studied in the past decades, and they continue to be a very active field of research. The interest in compressed sensing has motivated a significant part of the works on the topic. 
The starting point of many sparse reconstruction problems is a  \emph{dictionary}, 
intended as a collection of elements in the ambient space referred to as \emph{atoms} \cite{mallat}, used to represent the unknown quantity of interest.  The dictionary may be selected according to some \emph{a priori} information available on the problem of interest or, alternatively, its formation can be data-driven - see, e.g., \cite{dict1,dict2,dict3}. Recently, approaches aimed at learning the dictionary while jointly recovering the signal have also been developed \cite{chamb}. 
Typically, the cardinality of the dictionary is significantly larger than the dimension of the ambient space.
When the atoms in the dictionary do not form a basis for the ambient space, the dictionary is called \emph{redundant} or \emph{overcomplete}. The use of redundant dictionaries has proved to be a useful strategy in terms of artifact reduction, especially in the framework of signal denoising problems \cite{SED,SFM}.

Let $x\in\R^n$ be an unknown signal, and let $W=\{w_i\}_{i=1}^N$ be a dictionary with atoms $w_i\in\R^n$. We arrange the atoms as columns of the dictionary matrix, $\mW\in\R^{n \times N}$, with $n\ll N$, and refer to this matrix as the dictionary. In a synthesis perspective, a sparse reconstruction problem is the task of recovering a sparse vector $\alpha\in\R^N$, with most of its components vanishing,  that represents the original signal in terms of $\mW$, $x = \mW\alpha$, starting from a corrupted and possibly poorly sampled indirect observation  $b\in\R^{m}$ of $x$, with $m \leq n$. Assuming that the observation is linear in $x$ and the noise is additive, the sparse dictionary representation can be formulated as an optimization problem of the form
\begin{equation}
\mbox{minimize  $\|\alpha\|_0$  such that  $b = \mA x+\varepsilon$, with $x = \mW \alpha$,} 
\label{eq:cs}
\end{equation}
where $\|\alpha\|_0=\text{card}(\text{supp}(\alpha))$,  $\mA\in\R^{m \times n}$ is the forward model operator, and $\varepsilon\in\R^m$ is an additive noise vector. Addressing the minimization problem \eqref{eq:cs}  directly is a challenge due to its NP-hardness, thus explaining the need for alternative approaches. A strategy which has been widely explored and goes under the name of \emph{basis pursuit}  replaces the $\ell_0$-(semi)norm with its $\ell_1$ convex relaxation \cite{donoho}. 

When the signal itself is known to be compressible, i.e., $\mathrm{card}\left(\left\{x_i \mid |x_i| < \epsilon\right\} \right)\ll n$ with $\epsilon>0$ arbitrarily small, and the forward model operator $\mA$ satisfies the {\em restricted isometry property} (RIP) condition, an optimal bound for the error $\|x - \hat{x}\|_2$, with $\hat{x}$ denoting the recovered signal, has been derived. Moreover, if $x$ is sufficiently sparse, the signal can be recovered exactly \cite{candes1}. When the signal itself is not sparse, but it allows a sparse or compressible representation in a given dictionary, the exact recovery results still hold, provided that $\mA$ satisfies a {\em restricted isometry property adapted to a dictionary} (D-RIP) condition \cite{drip}. Given the theoretical motivation, a significant amount of research is devoted to identifying classes of operators to which these results could be applied. 

Besides the convex approaches, $\ell_p$-norms with $p<1$ have also been considered in place of the $\ell_0$-(semi)norm in problem \eqref{eq:cs}, as they are known to promote sparsity more strongly than the case $p=1$. Nonetheless, the presence of local minima is a clear limitation to the reliability of non-convex strategies. 

The sparse reconstruction problem for linear inverse problems allows a natural formulation in the Bayesian computational framework, with the notion of sparsity promoting priors.
In a number of previous works \cite{CSS,CPrSS,CPrS}, the  recovery of a sparse signal $x$ - or of a signal admitting a sparse representation in a given basis - has been addressed by modeling its entries in a hierarchical Bayesian framework as conditionally Gaussian random variables with unknown variances, with a generalized gamma hyperprior distribution.
The sparsity promotion and the convexity properties of the corresponding class of hypermodels have been  studied in \cite{CSS,CPrSS}. The derived results,
and in particular the considerations on the convexity properties of the resulting maximum a posteriori estimation problems, motivated the introduction of a hybrid hypermodel combining the strong sparsity promotion that typically characterizes non-convex settings with the convexity guarantees \cite{CPrS}.  


In this article, the hierarchical Bayesian framework outlined in the previous articles, and in particular the use of the hybrid algorithm of \cite{CPrS}, is extended to address sparse recovery problems in presence of redundant dictionaries. We consider a version of the iterated alternating sequential (IAS) algorithm that combines ideas from the Bayesian inference and iterative Krylov subspace methods, suitable for large scale problems, and is therefore particularly attractive for problems with large dictionaries. Numerical examples demonstrate the computational efficiency of the approach, and most importantly, show that for composite frame dictionaries, where each subdictionary would provide a sufficient representation of the signal, the method is capable of identifying an optimally sparse representation. 


\section{Hierarchical Bayesian formulation}
\label{sec:bayes}
Consider the linear inverse problem
\begin{equation}
\label{eq:model}
b = \mA x + \varepsilon\,,\; \varepsilon \sim \mathcal{N}(0,\mSigma) \quad\text{ such that }\quad x = \mW \alpha\,,
\end{equation}
where $\mA\in\R^{m \times n}$, with $m \leq n$, is the known forward model operator, $x\in\R^n$ is the unknown of interest and $\mSigma\in\R^{m \times m}$ is the symmetric positive definite covariance matrix of the additive Gaussian noise. In addition, we assume that $x$ admits a representation in the redundant dictionary $\mW\in \R^{n \times N}$, with $n \ll N$, where the unknown vector $\alpha\in\R^N$ is sparse, i.e. $\|\alpha\|_0 \ll N$, either because $x$ can be naturally described by few atoms in the dictionary or because $x$ needs to be compressed. 
\\In a number of previous contributions \cite{CSS,CPrSS,CPrS}, the \emph{a priori} sparsity belief on the unknown has been exploited by modeling its entries as independent random variables following a \emph{conditionally Gaussian} distribution, i.e.,
\begin{equation}
\label{eq:cond gauss}
\alpha_j\mid \theta_j \sim \mathcal{N}(0,\theta_j)\,,\; 1 \leq j \leq N\,,
\end{equation}
or, in an equivalent compact form,
\begin{equation}
\label{eq:cond gauss 2}
\alpha\mid \theta \sim \mathcal{N}(0,\mD_{\theta})\,,\; \mD_{\theta}=\mathrm{diag}\left(\theta_1,\ldots,\theta_N\right)\in\R^{N\times N}\,.
\end{equation}
The conditional Gaussian prior on $\alpha$ given the vector $\theta$ takes the form
\begin{align}
\label{eq:cond prior}
\pi_{\alpha\mid\theta} (\alpha\mid \theta) \:{\propto}\:&\displaystyle{ \frac{1}{\prod_{j=1}^{N}\sqrt{\theta_j}}\exp\bigg({-\frac{1}{2}\| \mD_{\theta}^{-1/2} \alpha \|^{2}}\bigg)}\\
\nonumber	
\:{=}\:&\displaystyle{\exp\bigg({-\frac{1}{2}\| \mD_{\theta}^{-1/2} \alpha \|^{2}} -\frac{1}{2}\sum_{j=1}^N \log\theta_j\bigg)}\,.
\end{align}
According to the Bayesian paradigm, the unknown vector of variances $\theta$ is also modeled as a random variable. The \emph{a priori} beliefs about $\theta$ are encoded in the \emph{hyperprior} $\pi_{\theta}(\theta)$, and the joint prior on the coupled vector of unknowns $(\alpha,\theta)$ reads
\begin{eqnarray}
\label{eq:joint pr}
\pi_{(\alpha,\theta)}(\alpha,\theta) = \pi_{\alpha\mid\theta}(\alpha\mid \theta)\pi_{\theta}(\theta)\,.
\end{eqnarray}
In \cite{CPrSS}, the authors propose to model the unknown variances $\theta_j$ as mutually independent random variables following a generalized gamma distribution,
\begin{equation}
\pi_{\theta}(\theta) = 
\pi_{\theta}(\theta \mid r,\beta,\vartheta) = \frac{|r|^n}{\Gamma(\beta)^n} \prod_{j=1}^{N}\frac 1{\vartheta_j} \left(\frac{\theta_j}{\vartheta_{j}}\right)^{r\beta - 1} \exp\bigg(-\left(\frac{{\theta_j}}{{\vartheta_{j}}}\right)^r\;\bigg)\,,
\label{eq:hyper}
\end{equation}
where $r\in \R\setminus\{0\}$, $\beta >0$, $\vartheta_j>0$. This choice is motivated by the observation that generalized gamma distributions tend to favor values which are close to the expected value while also allowing for few outliers very far from the mean. Presumably, the outlier variances give rise to the few non-zero values of $\alpha$, or values above a tiny threshold.

The information about the observation process is encoded in the likelihood distribution, which in view of the additive Gaussian noise model, takes the form
\begin{equation}
\label{eq:likelihood}
\pi_{b\mid \alpha}(b \mid \alpha) \propto \exp\left(-\frac{1}{2}\|\mS (\mA \mW \alpha - b)\|_2^{2}\right)\,,
\end{equation}
where $\mS$ is the Cholesky factor of the precision matrix $\mSigma^{-1}$, i.e. $\mSigma^{-1} = \mS^{\mT}\mS$. If the matrix $\mSigma$ and thereby $\mS$, are known, without loss of generality we can assume the noise to be white, i.e. $\mSigma = \mI$, because it can be whitened by a linear transform on $\mA$ and $b$, namely
\begin{equation}
\mA \longrightarrow \mS \mA\,,\quad b\longrightarrow \mS b\,.
\end{equation}
Under the white normal noise assumption, the likelihood distribution is of the form
\begin{equation}
\label{eq:likelihood 2}
\pi_{b\mid \alpha}(b \mid \alpha) \propto \exp\left(-\frac{1}{2}\| \mA \mW \alpha - b\|_2^{2}\right)\,.
\end{equation}
The conditional prior and the hyperprior are coupled to the posterior distribution via Bayes' formula, yielding the following expression for the posterior distribution
\begin{equation}
\label{eq:bayes}
\pi_{(\alpha,\theta)\mid b}(\alpha,\theta \mid b) \propto \pi_{b \mid \alpha}(b \mid \alpha)\pi_{(\alpha,\theta)}(\alpha \mid \theta) \pi_{\theta}(\theta).
\end{equation}
In the Bayesian framework, the posterior distribution is the complete solution to the inverse problem, that can be used to produce representative estimates of the unknown of interest, and quantify the uncertainty.
Here, we chose to summarize the posterior with the Maximum A Posteriori (MAP) estimate,
\begin{equation}
(\alpha^{*},\theta^*) \in \arg\max_{\alpha,\theta}\left\{\pi_{(\alpha,\theta)\mid b}(\alpha,\theta \mid b)\right\}
\end{equation}
or equivalently, by taking the negative logarithm of the density and ignoring the additive constants,
\begin{equation}
\label{eq:map}
(\alpha^{*},\theta^*) \in \arg\min_{\alpha,\theta}\left\{\mathcal{F}(\alpha,\theta)\right\},
\end{equation}
where
\begin{align}
\nonumber
&\mathcal{F}(\alpha,\theta) = \mathcal{F}(\alpha,\theta \mid r,\vartheta,\beta)\\
=&\quad \lefteqn{ \phantom{\frac 12 \|b - \mA  \mW\alpha\| ^2 + \frac 12 \sum_{j=1}^N\frac{\alpha_j^2}{\theta_j} }}
{\displaystyle\frac 12} \|b - \mA  \mW\alpha\|^2+
\underbrace{\frac 12\sum_{j=1}^N\frac{\alpha_j^2}{\theta_j}
	-\eta
	\sum_{j=1}^N\log\frac{\theta_j}{\vartheta_j} +\sum_{j=1}^N\left(\frac{\theta_j}{\vartheta_j}\right)^r}_{\mathcal{P}(\alpha,\theta \mid r,\beta,\vartheta)}\,,\quad \eta = \bigg(r\beta - \frac 32\bigg). \label{eq:obj map}
\end{align}
Echoing the terminology of classical regularization schemes, we refer to $\mathcal{P}(\alpha,\theta \mid r,\beta,\vartheta)$ as the \emph{penalty term}.

\section{The IAS algorithm}
\label{sec:ias}
The search for the minimizer of the MAP objective function in \eqref{eq:obj map} is carried out with the {\em global hybrid} scheme introduced in \cite{CPrS}, based  on the iterative alternating sequential (IAS) algorithm described below.  Details of the hybrid scheme that ensues are reviewed in Section \ref{sec:ghias}.

Given a suitable initialization of the variances $\theta^{0}$,  at each iteration step the IAS algorithm updates the iterates $\alpha^t$, $\theta^t$ by solving the minimization problem in alternating directions, that is
\begin{equation}
\label{eq:ias}
\alpha^{t+1}\in\arg\min_{\alpha}\left\{\mathcal{F}(\alpha,\theta^{t})\right\}\,,\quad 
\theta^{t+1}\in\arg\min_{\theta}\left\{\mathcal{F}(\alpha^{t+1},\theta)\right\}\,.
\end{equation}
Because of the particular form of the objective function, both variables can be updated efficiently as follows.

\paragraph{\textbf{Update of }$\bm{\alpha}$}
The $\alpha$-update reduces to solving a quadratic minimization problem, i.e.,
\begin{equation}
\label{eq:get alpha}
\alpha^{t+1}\in\arg\min_{\alpha}   \left\{ \|b-\mA\mW\alpha\|_2^2  + \|\mD_\theta^{-1/2}\alpha\|^2 \right\},\quad \theta = \theta^t,
\end{equation}
or, equivalently, finding the solution in the least squares sense of the linear system
\begin{equation}
\label{eq:get alpha LS}
\begin{bmatrix}
\mA\mW\\
\mD_{\theta}^{-1/2}
\end{bmatrix}\alpha = \begin{bmatrix}
b\\
0
\end{bmatrix}\,.
\end{equation}
After performing the change of variable
\begin{equation}
\mD_{\theta}^{-1/2}\alpha = \gamma\,,
\end{equation}
we can write  \eqref{eq:get alpha LS} as
\begin{equation}
\label{eq:get gamma}
\begin{bmatrix}
\mA\mW\mD_{\theta}^{1/2}\\
\mI
\end{bmatrix}\gamma = \begin{bmatrix}
b\\
0
\end{bmatrix}\,,
\end{equation}
where $\mI$ is an $n\times n$ unit matrix. The solution of this least squares problem is also Tikhonov regularized solution of 
\begin{equation}
\mA\mW\mD_{\theta}^{1/2}\gamma = b,\,\; \alpha = \mD_{\theta}^{1/2}\gamma\,,
\end{equation}
with regularization parameter equal to one. An alternative to Tikh\-on\-ov regularization yielding a similar solution is to solve the underlying linear system with an iterative solver equipped with an early stopping criterion. The stopping condition is usually based on a variant of Morozov discrepancy principle, whereas the iterations terminate as soon as the discrepancy is of the order of the observation noise. In the statistical framework, under the Gaussian noise assumption, the noise level can be expressed in terms of the standard deviation of the noise.  In our case, where we assume  $m$-dimensional white noise, this quantity is equal to $\sqrt{m}$.   Following \cite{SIREV,CPrSS}, we solve the linear system using the Conjugate Gradient for Least Squares (CGLS) algorithm with the early stopping at noise level $\sqrt{m}$; see \cite{SIREV}  for more details.

\paragraph{\textbf{Update of }$\bm{\theta}$} Due to the mutual independence of the entries of $\theta$, each variance $\theta_j$ can be updated separately by imposing the 
component-wise first order optimality condition on \eqref{eq:obj map}. More specifically, $\theta_j^{t+1}$ is the solution of the non-linear equation
\begin{equation}
\label{eq:get th}
\frac{\partial{\mathcal F}}{\partial\theta_j} = 
-\frac{1}{2} \frac{\alpha_j^2}{\theta_j^2}-\left(r \beta -\frac{3}{2}\right)\frac{1}{\theta_j} + r \frac{\theta_j^{r-1}}{\vartheta_j^r} = 0\,,\; \alpha=\alpha^{t+1}\,.
\end{equation}
For some values of $r$, e.g., $r=\pm 1$, \eqref{eq:get th} admits an analytic solution. However, in general we need to solve it numerically. It was shown in \cite{CPrSS} that after the changes of variables  $\theta_j = \vartheta_j \xi_j$, $\alpha_j = \sqrt{\vartheta_j} z_j$, we may write $\xi_j = \varphi(|z_j|)$, and  via implicit differentiation, the function $\varphi$ satisfies the initial value problem
\begin{equation}\label{get xi 2}
\varphi'(z) = \frac{2z \varphi(z)}{2r^2 \varphi(z)^{r+1}+z^2} , \quad \varphi(0) = \left(\frac{\eta}{r}\right)^{1/r}.
\end{equation}
Therefore the updated value of $\theta_j$ can be computed by a numerical time integrator. Since the same type of differential equation is satisfied by all components, an efficient way to update $\theta$ is to sort the current values $z_j$ in an ascending order, and integrate sequentially over the gaps between the values by a suitable time integrator.

We point out that unlike in the formally similar alternating direction method for multipliers (ADMM) algorithm \cite{boyd2011distributed} that is often used to solve regularized inverse problems with sparsity promoting priors, the IAS algorithm does not require the introduction of an artificial decoupling term of the fidelity and penalty terms, as the partial decoupling in IAS is automatic and exact.

From the point of view of statistical analysis, the proposed algorithm is not aiming at exploring the posterior density, and the MAP estimate might not be the best single point estimate to characterize the posterior. Other strategies of interest include the marginalization of the posterior density with respect to the hyperparameter $\theta$, or estimating an optimal $\theta$ by first marginalizing $\alpha$. These alternative strategies have been discussed in literature, see, e.g., \cite{Vidal} for a recent and comprehensive contribution.

\section{Parameter selection strategies}
Before presenting the details of the hybrid scheme used in the numerical tests, we briefly review some of the main results related to the selection  of the hyperparameters $(r,\beta,\vartheta)$ appearing in the expression of the hyperprior in \eqref{eq:hyper}. 

We start recalling a theorem, whose proof can be found in \cite{CPrSS}, summarizing how $r$ and $\beta$ affect the convexity properties of the functional $\mathcal{F}$ .

\begin{theorem}\label{th:convexity}
	Let $\beta>0$ and $r\neq 0$, and let $\mathcal{F}(\alpha,\theta)$ be the objective function for the minimization problem in (\ref{eq:map}).
	\begin{itemize}
		\item[(a)] 
		If $r\geq 1$ and $\eta = r\beta  -3/2>0$, the function $\mathcal{F}(\alpha,\theta)$ is globally convex.
		\item[(b)] 
		If $0<r<1$ and $\eta = r\beta  -3/2>0$, or, if $r<0$ and $\beta>0$,
		the function $\mathcal{F}(\alpha,\theta)$ is convex provided that
		\begin{equation} \label{th bar}
		\theta_j < \overline{\theta}=\vartheta_j \left(\frac{\eta}{r|r-1|}\right)^{1/r}.
		\end{equation}
	\end{itemize}
\end{theorem}

The convexity of the MAP objective function, guaranteed for $r \geq 1$, is very convenient, however some of the configurations attained for $r<1$ can be very attractive in terms of sparsity promotion and rate of convergence. To better understand the connection between the parameter of the hyperprior and sparsity, consider the updating formula  \eqref{get xi 2}, expressing $\theta_j$ as a function of $\alpha_j$ as
\begin{equation}\label{def g}
\theta_j =  g_j(\alpha_j) =  \vartheta_j \varphi\left(\frac{|\alpha_j|}{\sqrt{\vartheta_j}}\right)\,.
\end{equation}
We review some recent results \cite{CHPS,CSS,CPrSS} about the connections that can be drawn between the generalized gamma hyperpriors and classical sparsity promoting penalty terms, assuming that $(\theta_j,\alpha_j)$ satisfies the above identity.
\begin{itemize}
	\item[(i)] For the gamma hypermodel, i.e. $r=1$, as $\eta=r\beta-\frac{3}{2}\to 0^{+}$ the penalty term approaches a weighted $\ell_1$-penalty term \cite{CHPS,CSS},
	\begin{equation}
	\label{eq:lim ell1}
	\lim_{\eta \to 0^+}\mathcal{P}\left(\alpha,g(\alpha)\mid 1,\frac{3}{2}+\eta,\vartheta\right) = \sqrt{2} \sum_{j=1}^n \frac{|\alpha_j|}{\sqrt{\vartheta_j}}\,.
	\end{equation}
	\item[(ii)] If $r\beta=\frac{3}{2}$, 
	the penalty term coincides with the weighted $\ell_p$-norm, with $p =2r/(r+1)$ \cite{CPrSS},
	\begin{equation}
	{\mathcal P}\left(\alpha, g(\alpha)\mid r,\frac{3}{2r},\vartheta\right)  = C_r\sum_{j=1}^n \frac{|\alpha_j|^p}{\sqrt{\vartheta_j}^p}, \quad C_r = \frac{r+1}{(2r)^{r/(r+1)}}.
	\end{equation}
	\item[(iii)]For the inverse gamma hypermodel, corresponding to $r=-1$, the penalty term approaches the Student distribution, a prominently fat tailed distribution favoring large outliers, and leading to a greedy algorithm that strongly promotes sparsity \cite{CPrSS}. 
\end{itemize}

To summarize, the above results indicate that
the hyperpriors for which the global convexity of the corresponding hypermodel is not guaranteed ($r<1$) are expected to promote sparsity more effectively than the limit case $r=1$ that can be seen as a counterpart of the $\ell_1$-penalized case.

While the hyperparameters $r$ and $\beta$ determine the strength of the sparsity promotion and the convexity properties of the MAP objective function, the vector of the scale parameters $\vartheta$ can be set automatically once the operator $\mA\mW$ is given. More specifically, for each $j$, $\vartheta_{j}$ can be related to the sensitivity of the data to $x_j$, given by the quantity  $\|\mA\mW e_j\|^2$, where $e_j\in\R^N$ denotes the canonical $j$-th Cartesian unit vector. It was proven, for $r=1$ in \cite{SIREV,CSS} and in more general settings in \cite{CPrSS}, that under the assumption that the signal-to-noise ratio is given, and that the prior satisfies an exchangeability condition guaranteeing that no particular sparse combinations of components of $x$ are favored over others, 
the entries of $\vartheta$ must be chosen as
\begin{equation}
\label{eq:sens}
\vartheta_j = \frac{C}{\|\mA\mW e_j\|^2}\,,
\end{equation}
where $C>0$ is a constant encoding the expected sparsity on the solution and an estimate of the signal-to-noise-ratio. For details, we refer to the cited articles.
We remark that, in general, sensitivity weights are introduced to compensate for the possible non-uniform design of the forward model operator $\mA$.  Sensitivity weights play an important role in, e.g.,  inverse source problems, in which sources near the observation points may be favored over far away sources unless the exchangeability condition is imposed. In the current setting, when the dictionary consists of sub-frames with possibly different column norms, we expect the different weights $\vartheta_j$ to prevent the representation of the signal 
the frames with larger column norms to dominate.

\section{Local and global hybrid IAS}
\label{sec:ghias}

In the following discussion, we write the penalty function ${\mathcal P}(x,\theta\mid r,\beta,\vartheta)$ in terms of components,
\begin{eqnarray*}
	{\mathcal P}(\alpha,\theta\mid r,\beta,\vartheta) &=& \sum_{j=1}^N\left(\frac 12 \frac{\alpha_j^2}{\theta_j}  -\eta\log\frac{\theta_j}{\vartheta_j} +\left(\frac{\theta_j}{\vartheta_j}\right)^r\right)\\
	&=&\sum_{j=1}^N {\mathcal P}_j(\alpha_j,\theta_j\mid r,\beta,\vartheta_j).
\end{eqnarray*}

In \cite{CPrS}, two different hybrid strategies to speed up and enhance sparsity promotion in the IAS algorithm were proposed. In both versions, the IAS iterations are initiated by selecting a conservative set of hyperparameters for which the objective function is convex, thus guaranteeing global convergence to a unique minimizer. We denote this set of parameters by $(r^{(1)},\beta^{(1)},\vartheta^{(1)})$. For the second phase of the hybrid algorithm, we select another set of parameters, $(r^{(2)},\beta^{(2)},\vartheta^{(2)})$, for which the global convexity of the objective function is not valid.
To match the models so that they express coherent prior beliefs, we adjust the scale parameters $\vartheta^{(j)}$ so as to satisfy the compatibility condition
\begin{equation}\label{compatibility}
\left(\frac{\eta^{(1)}}{r^{(1)}}\right)^{1/r^{(1)}}\vartheta_j^{(1)} = \left(\frac{\eta^{(2)}}{r^{(2)}}\right)^{1/r^{(2)}}\vartheta_j^{(2)},
\end{equation}
that guarantees that the parameter $\theta_j$ computed at $\alpha_j=0$ returns the same value regardless of the model. 
For further discussion, we refer to \cite{CPrS}.

In the local hybrid version,  the IAS algorithm is initially run with hyperparameters  $(r^{(1)},\beta^{(1)},\vartheta^{(1)})$, and after each iteration step, we check which $\theta_j$, if any, satisfies the condition (\ref{th bar}), where $\overline\theta$ is computed using the hyperparameter set $(r^{(2)},\beta^{(2)},\vartheta^{(2)})$. In correspondence of those which do, we modify the local objective function so that 
\begin{equation}
{\mathcal P}_j(\alpha_j,\theta_j\mid r^{(1)},\beta^{(1)},\vartheta_j^{(1)}) \rightarrow {\mathcal P}_j(\alpha_j,\theta_j\mid r^{(2)},\beta^{(2)},\vartheta_j^{(2)}).
\end{equation}

The global hybrid scheme is based on the idea that after a number of IAS iteration rounds, the iterate of the globally convex objective function with hyperparameters
$(r^{(1)},\beta^{(1)},\vartheta^{(1)})$  is near the unique global minimum of that objective function. Restarting the IAS from the current point with the parameters $(r^{(2)},\beta^{(2)},\vartheta^{(2)})$ may quickly find a local minimizer near the global minimizer of the original objective function. While the two minimizers are likely not far apart, the local minimizer is typically sparser, and the convergence to it is faster.


\section{Computed Examples}

In this section, we demonstrate the viability of the hybrid IAS algorithm in the context of overcomplete representations. More specifically,  we restrict ourselves to the global hybrid strategy, switching from the first to the second hyperprior after 10 iterations, if not differently specified. The main goal of the following examples is to demonstrate that the global hybrid IAS is capable of selecting from a dictionary of sub-frames, where several representations are admissible, a set of atoms that make the representation as sparse as possible.

\medskip

\begin{figure}[!t]
	\centering
	\centerline{		\includegraphics[scale = 0.25]{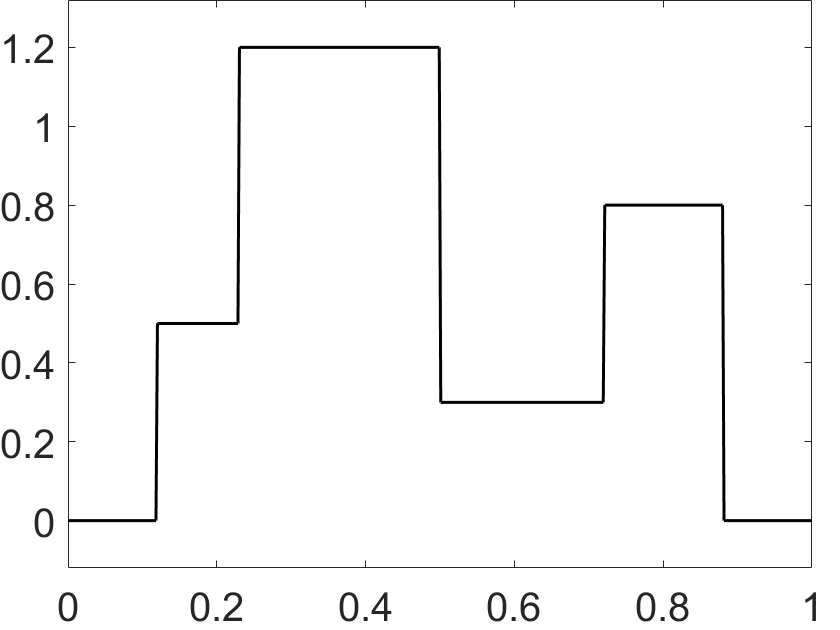}\quad
		\includegraphics[scale = 0.25]{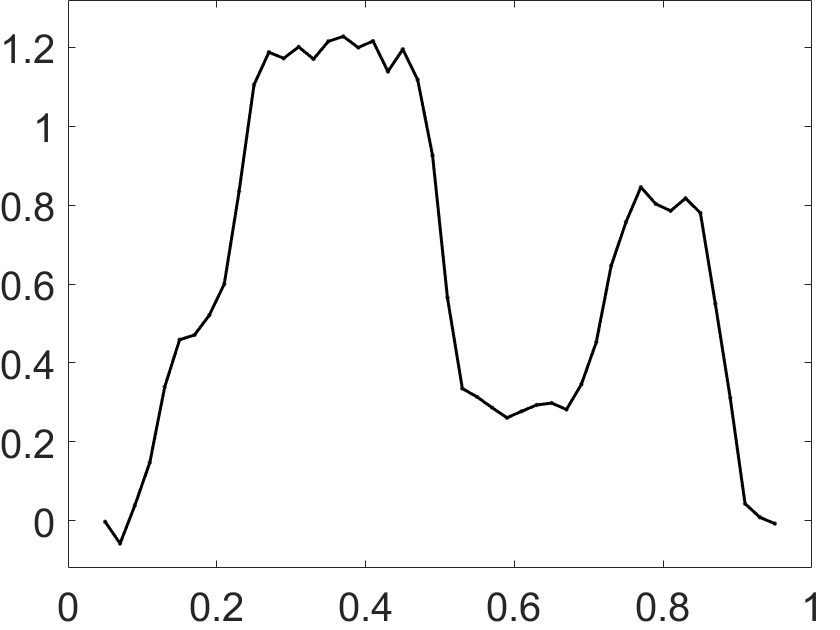}	}
	\caption{The generative model (left) and the blurred and noisy data vector $b \in \R^{46}$ (right).}
	\label{fig:1D_Data}
\end{figure}

\paragraph{\textbf{Signal restoration from convolution data}}
The first test case is a one-dimensional deconvolution problem. The generative model is a piecewise constant signal  $f:[0,1]\to \R$,  $f(0)=0$, and the data consist of a few discrete observations, 
\begin{equation}
b_j = \int_0^1 A(s_j - t) f(t) dt + \varepsilon_j,\quad 1\leq j\leq m,  \quad A(s_j,t) = \frac{1}{\sqrt{2\pi w^{2}}}e^{-\frac{(s_j - t)^{2}}{2w^{2}}}, 
\end{equation}
corrupted by Gaussian blur with $w=0.02$ and additive scaled white Gaussian noise, with standard deviation $\sigma$ set to $2\%$ of the maximum of the noiseless signal.
The data has been generated using a discretization of the unit interval with $n=n_{\rm dense} = 1253$ nodes, while in the forward model used for solving the inverse problem, we set $n=500$. The number of equidistant observation points in the signal domain is $m=46$. The generative signal and the data are shown in Figure~\ref{fig:1D_Data}.

The generative signal admits a natural sparse representation in terms of its increments  $z_j = z_j - z_{j-1}$ over the interval of definition. Assuming  $x_0=0$, then 
\begin{equation}\label{difference}
z = \mB  x\,,\quad \mB=\left[\begin{array}{cccc}
1&0&\ldots&0\\
-1&1&\ldots&0\\
&&\ddots&\\
0&\ldots&-1&1\\
\end{array}\right]
\in \R^{n\times n },
\end{equation}
hence
\begin{equation}
x = \mL z\quad\mathrm{with}\quad  \mL = \mB^{-1} = \left[\begin{array}{cccc}
1&0&\ldots&0\\
1&1&\ldots&0\\
\vdots &&\ddots&\\
1&\ldots&1&1\\
\end{array}\right]\in\R^{n\times n}.
\end{equation}
Our goal is to test the effectiveness of the outlined framework in recovering the most natural sparse representation of the given signal. 
Let $\mC$ denote the discrete cosine transform matrix, providing an alternative and accurate way of representing the signal,
\begin{equation}
x =  \mC^\mT y, \quad y = \mC x,
\end{equation}
which is, however, not sparse. To test wether the algorithm is able to identify the frame that allows a sparse representation,
we consider the overcomplete dictionary,
\begin{equation}
\label{eq:dict 1d}
\mW= [\mW_1\,,\,\mW_2]\in\R^{n \times 2n} \quad\text{with}\quad \mW_1 = \mL \in \R^{n \times n }\;\text{and}\;\mW_2 = \mC^{\mT}\in\R^{n \times n}\,,
\end{equation}
and formulate the underlying linear inverse problem as
\begin{equation}
b = \mA\mW \alpha + \varepsilon = \mA [\mW_1\,,\,\mW_2] \begin{bmatrix}
\alpha_1\\
\alpha_2
\end{bmatrix}+\varepsilon, \quad \varepsilon \sim \mathcal{N}(0,\sigma^2 \mI_{m})\,,
\end{equation}
where $\mA$ is the discrete blur operator.

In this example, the global hybrid IAS is run with parameters $(r^{(1)},\eta^{(1)})=(1,10^{-4})$, $(r^{(2)},\eta^{(2)})=(1/2,10^{-3})$; we recall that the sensitivity weights $\vartheta^{(1)}$ are set automatically according to \eqref{eq:sens}, while the vector $\vartheta^{(2)}$ is fixed so that condition \eqref{compatibility} is satisfied.

\begin{figure}[!t]
	\centering
	\centerline{
		\includegraphics[scale = 0.23]{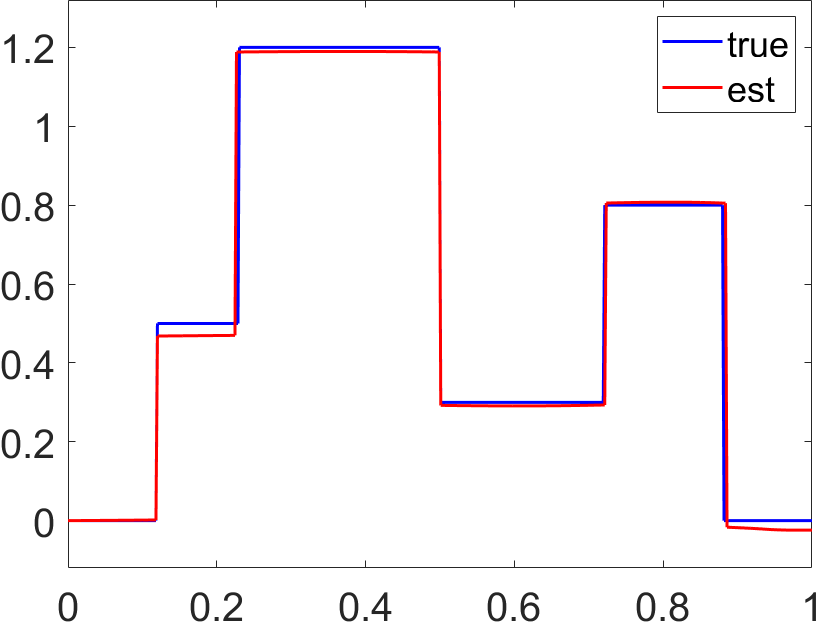}\quad 
		\includegraphics[scale = 0.23]{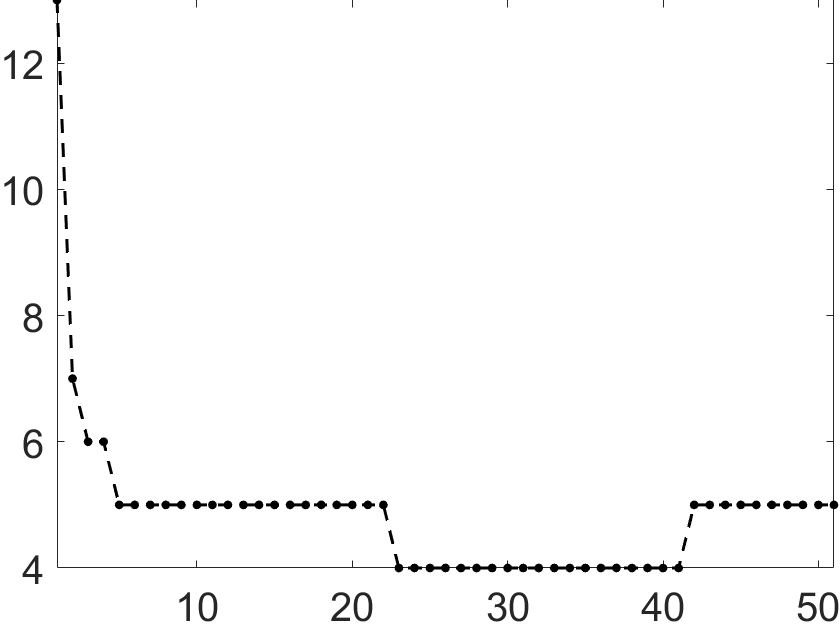}
	}
	\caption{Reconstruction of the signal $x$ (left) and the count of CGLS steps per outer iteration of the global hybrid IAS (right).}
	\label{fig:rec 1d}
\end{figure}

\begin{figure}[!t]
	\centering
	\centerline{\includegraphics[width = 4.3cm]{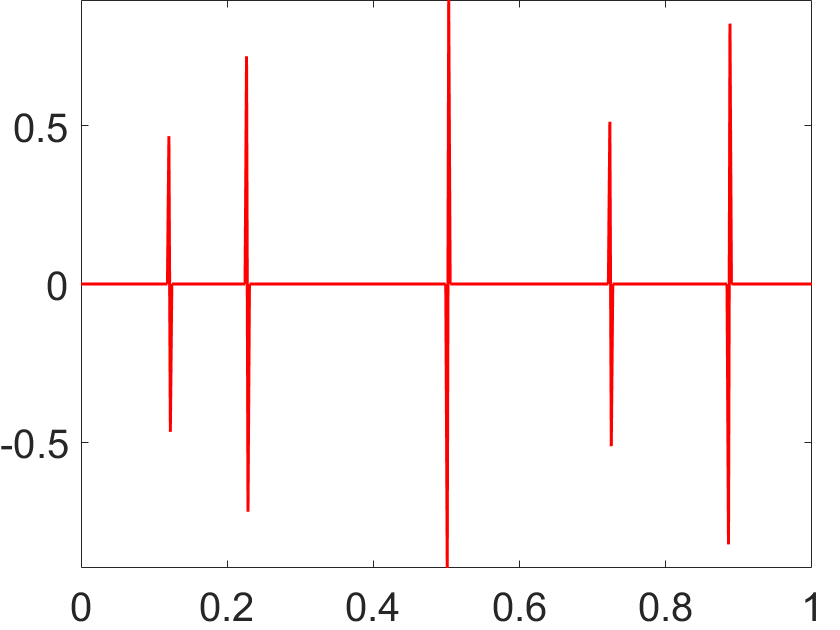}		\includegraphics[width = 4.3cm]{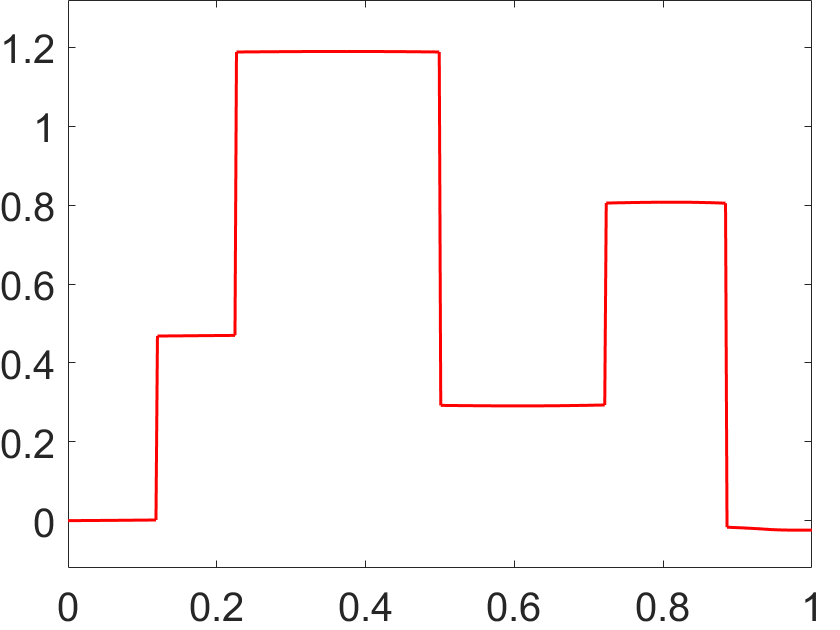}		\includegraphics[width = 4.3cm]{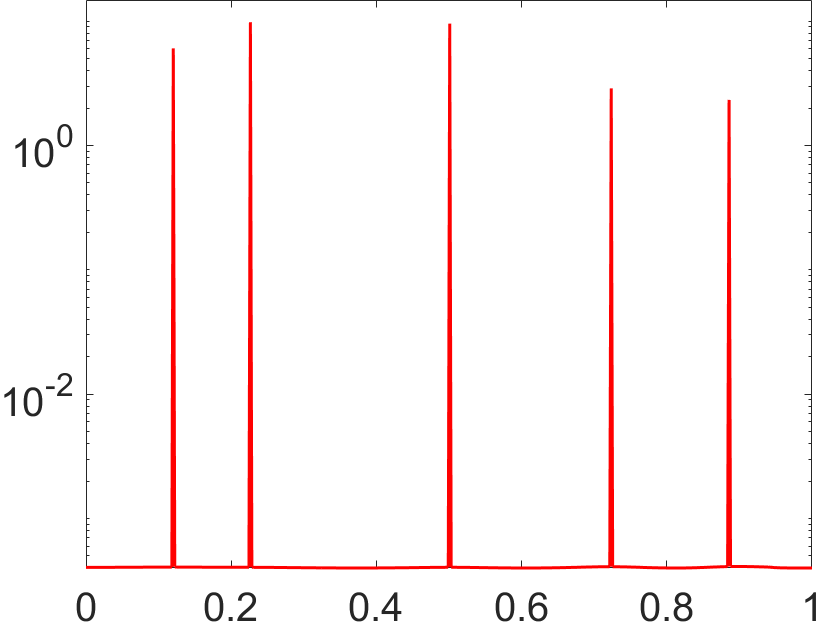}
	}
	\centerline{	\includegraphics[width=4.3cm]{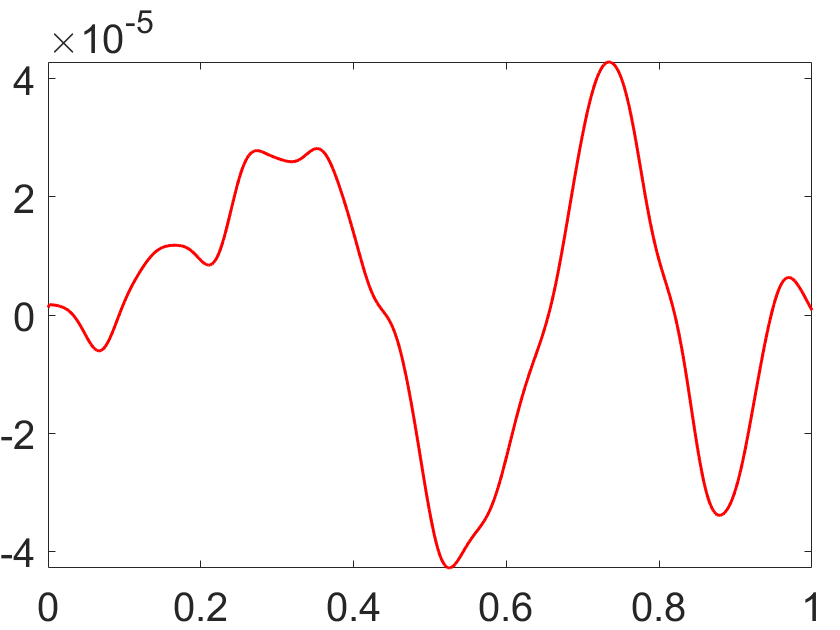}	\includegraphics[width=4.3cm]{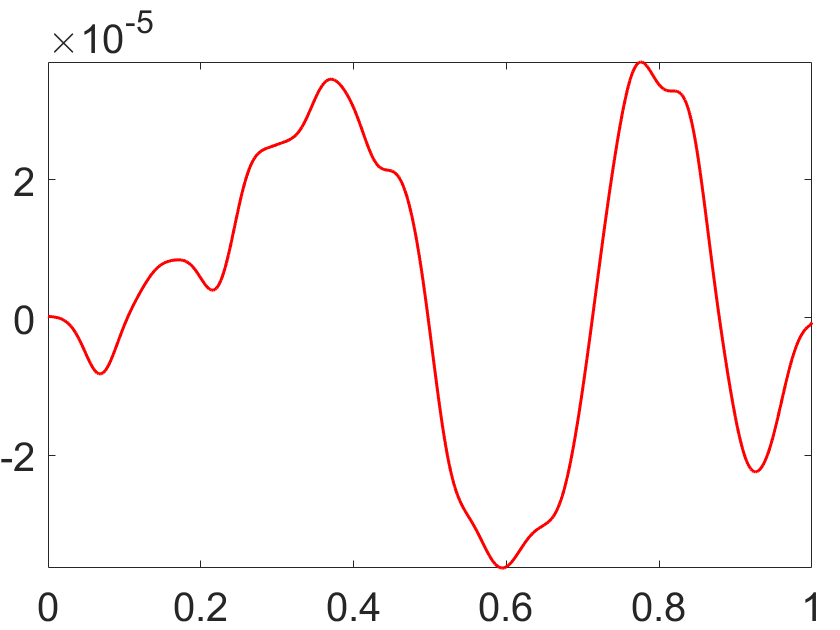}		\includegraphics[width=4.3cm]{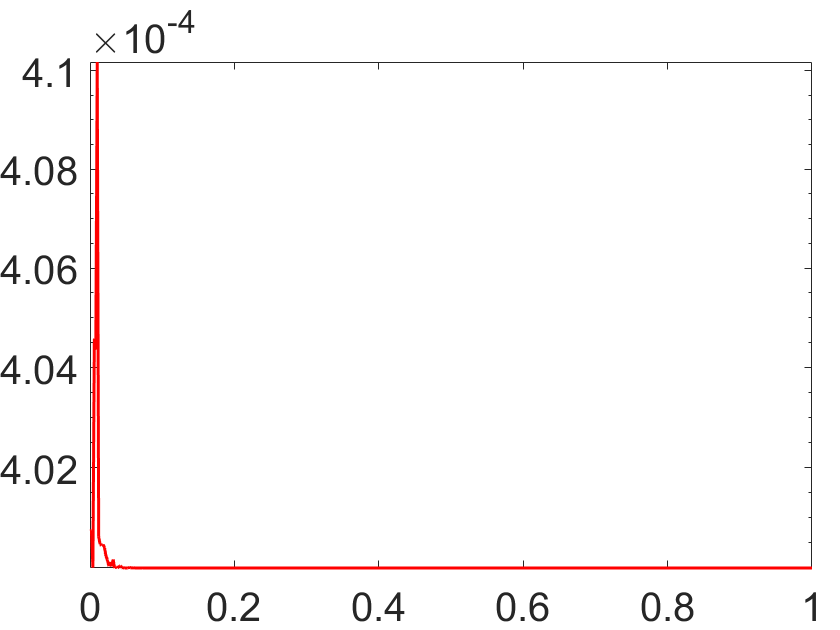}	
	}
	\caption{Vector $\alpha_i$ (left panles), corresponding scaled variances (middle panels) and contribution of the signal $\mW_i \alpha_i$ (right panels) for $i=1$, i.e. representation in terms of increments, (top row) and $i=2$, i.e. representation in terms of cosine transform, (bottom row).}
	\label{fig:rec 1d 2}
\end{figure}

The signal reconstructed by the global hybrid IAS scheme is shown in Figure \ref{fig:rec 1d}. The restored $\alpha_1$ and $\alpha_2$ and their contribution in the estimated signal are shown in Figure \ref{fig:rec 1d 2}, together with the scaled variances corresponding to $\alpha_1$ and $\alpha_2$, i.e.
\begin{equation}
\frac{\theta_j}{\vartheta_j^{(2)}}\,,\;1 \leq j \leq n \quad\text{and}\quad \frac{\theta_j}{\vartheta_j^{(2)}}\,,\;n+1 \leq j \leq 2n\,.
\end{equation}
Notice that the output variances are scaled by the sensitivities corresponding to the second hyperprior  used to design the hybrid scheme.

\noindent Despite the relatively high level of degradation (blur and noise) and down-sampling in the observed data $b$, the algorithm has no problem detecting the basis that provides a more natural and sparse representation for the original signal. In fact, the coefficients $\alpha_2$ are five to six orders of magnitude smaller than the non-vanishing components of $\alpha_1$. The degree of sparsity in the final representation is also reflected in the number of CGLS steps per outer iteration of the global hybrid IAS - see Figure \ref{fig:rec 1d} - which quickly settles around the cardinality of the support of $\alpha$. 

\medskip

\begin{figure}[!t]
	\centering
	\centerline{
		\includegraphics[height=3.75cm]{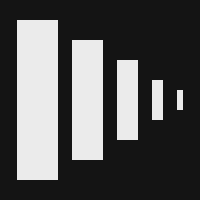}\;
		\includegraphics[height=3.75cm]{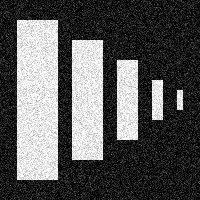}\;	\includegraphics[height=3.75cm]{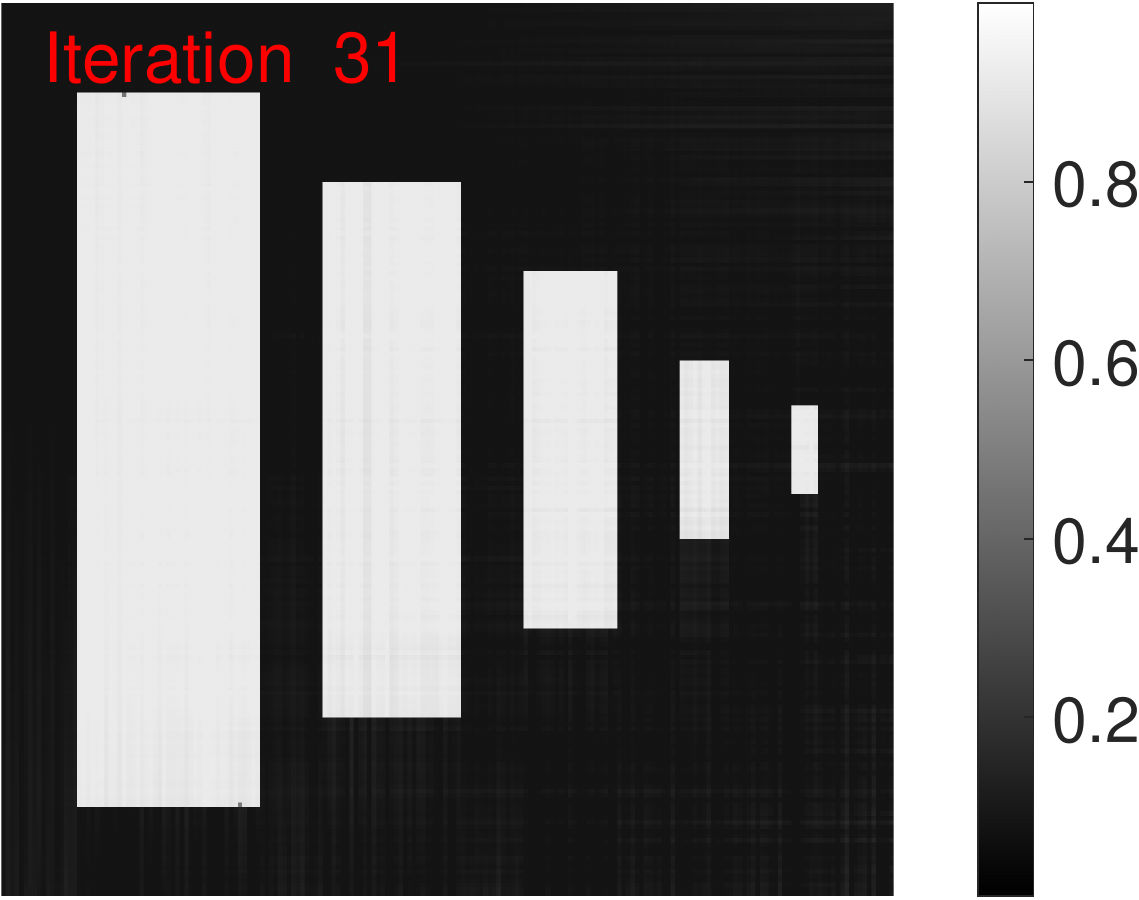}	}
	\caption{Original image (left), observed data (middle) and reconstructed image (right).}5
	\label{fig:Data_2d}
\end{figure}

\paragraph{\textbf{Image denoising on a synthetic image}}
In the second example, we consider the problem of denoising a blocky gray scale test image $x\in\R^{n \times n}$, $n=200$. The pixel values, which are between $0$ and $1$, are corrupted by scaled white Gaussian noise with standard deviation $\sigma$ set to $10\%$ of the maximum of the noiseless image, i.e. $\sigma=0.1$ - see Figure \ref{fig:Data_2d}.

The test image presents sharp edges lying along the horizontal and vertical axes. Therefore, $x$ admits a sparse representation both in the vertical and horizontal increment bases, the latter being slightly less sparse than the former. After representing the image in vector form $x\in\R^{n^2}$ by stacking the pixel values columnwise, we introduce the redundant dictionary $\mW= [\mW_1\,,\,\mW_2]\in\R^{n^2 \times 2n^2}$ with
\begin{equation}
\label{eq:dict 2d}
\mW_1 = \left(\mI_{n} \otimes \mB\right)^{-1} \in \R^{n^2 \times n^2 }\;\text{and}\;\mW_2 = \left(\mB \otimes \mI_n\right)^{-1}\in\R^{n^2 \times n^2}\,,
\end{equation}
where $\mB$ is defined as in \eqref{difference}, and $\otimes$ stands for the Kronecker product.  Homogenous Dirichlet boundary conditions are assumed on the left and top edges of the image. We want to estimate the sparse vector $\alpha = [\alpha_1\,,\,\alpha_2]^{\mT}$, with $\alpha_i\in\R^{n^2}$, $i=1,2$, from the data vector $b\in\R^{n^2}$, given the forward model
\begin{equation}
b = \mW \alpha + \varepsilon = [\mW_1\,,\,\mW_2] \begin{bmatrix}
\alpha_1\\
\alpha_2
\end{bmatrix}+\varepsilon, \quad \varepsilon \sim \mathcal{N}(0,\sigma^2 \mI_{n^2})\,.
\end{equation}
It is worth remarking here that we require $\alpha$ to be not only sparse, but as sparse as possible.

The hyperparameters of the global hybrid IAS are set as $(r^{(1)},\eta^{(1)})=(1,10^{-3})$ and $(r^{(1)},\eta^{(1)})=(1/2,10^{-2})$, while, as before, $\vartheta^{(1)},\vartheta^{(2)}$ are automatically fixed according to \eqref{eq:sens} and \eqref{compatibility}, respectively.

\begin{figure}[!t]
	\centering
	\centerline{\includegraphics[height=3.75cm]{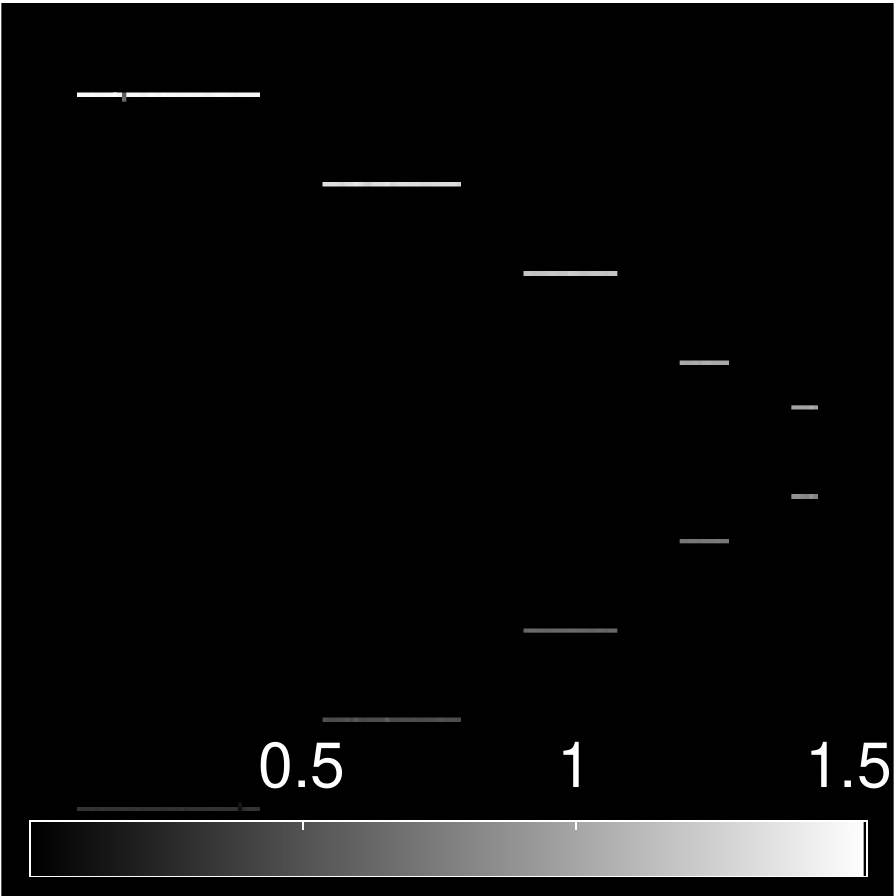}\;		\includegraphics[height=3.75cm]{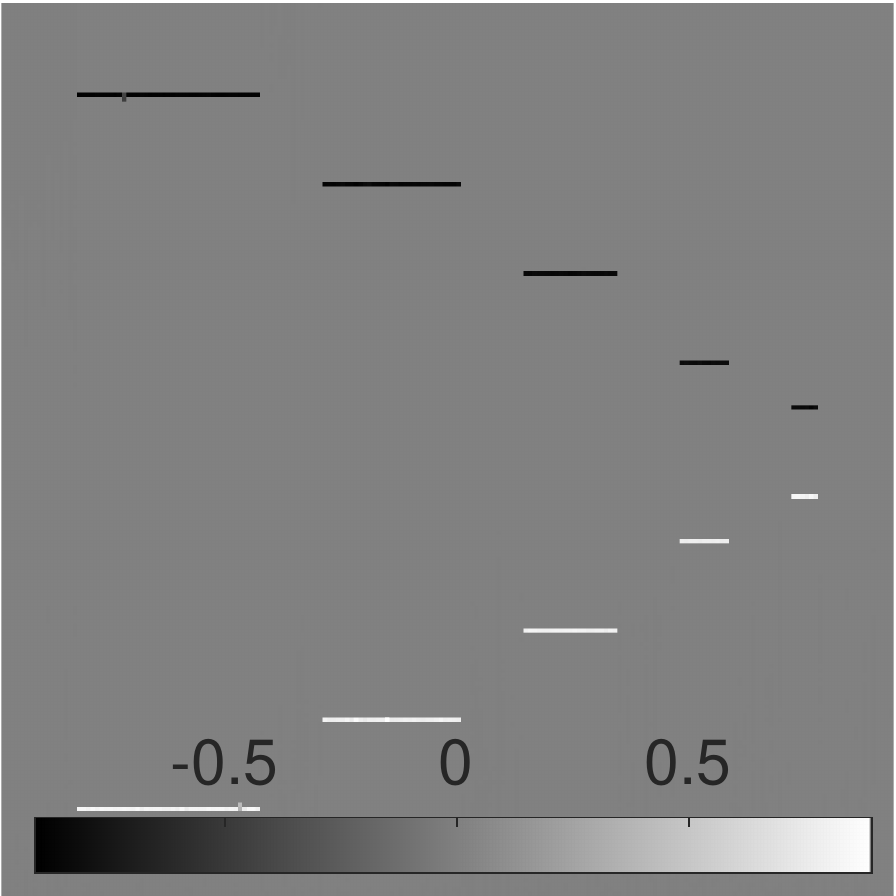}\; 		\includegraphics[height=3.75cm]{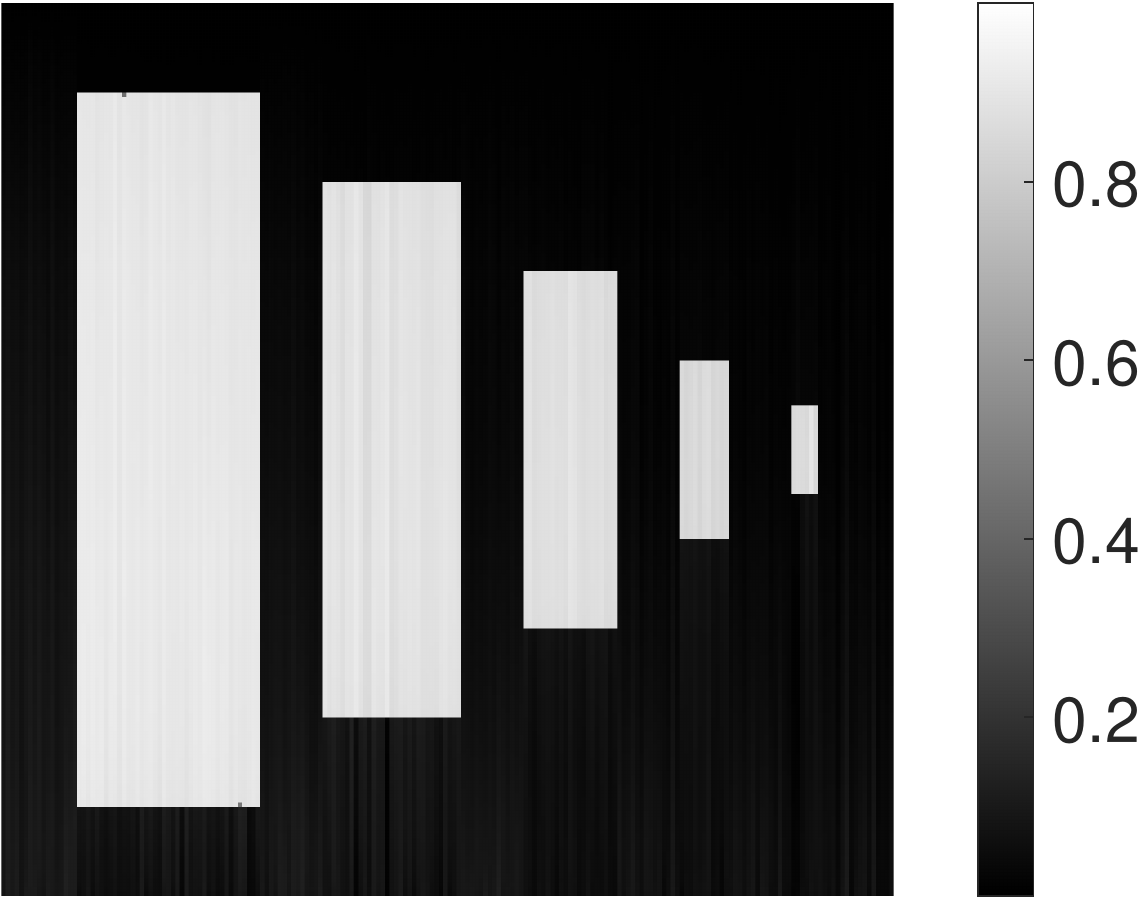} }
	\centerline{\includegraphics[height=3.75cm]{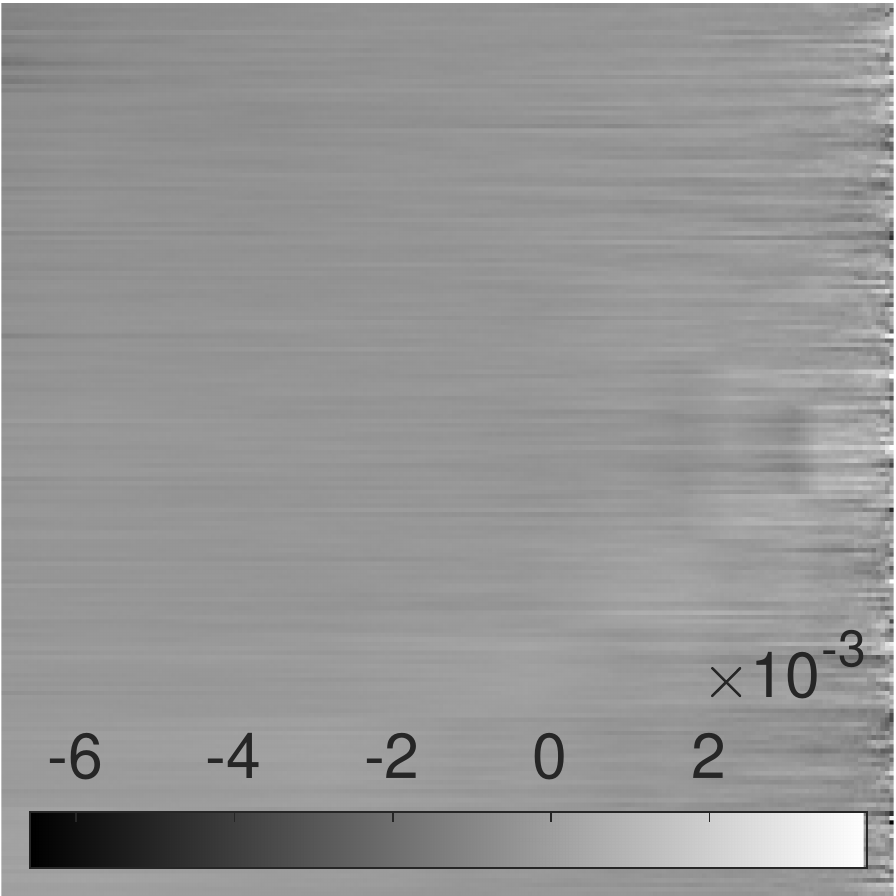}\;	\includegraphics[height=3.75cm]{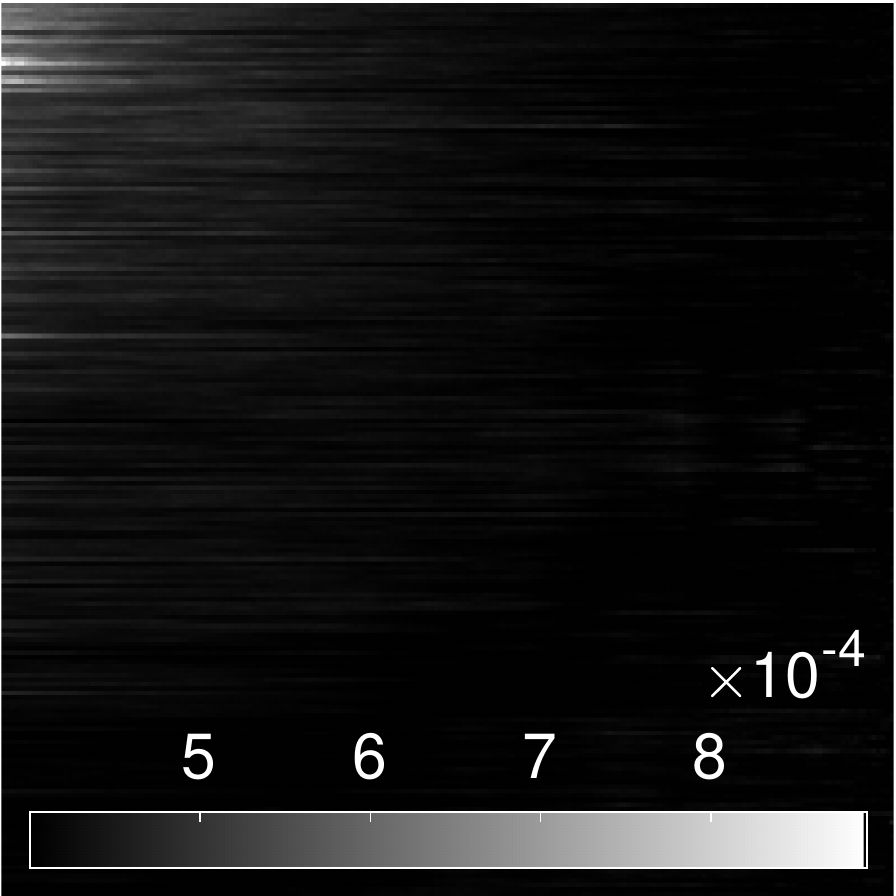}\; 		\includegraphics[height=3.75cm]{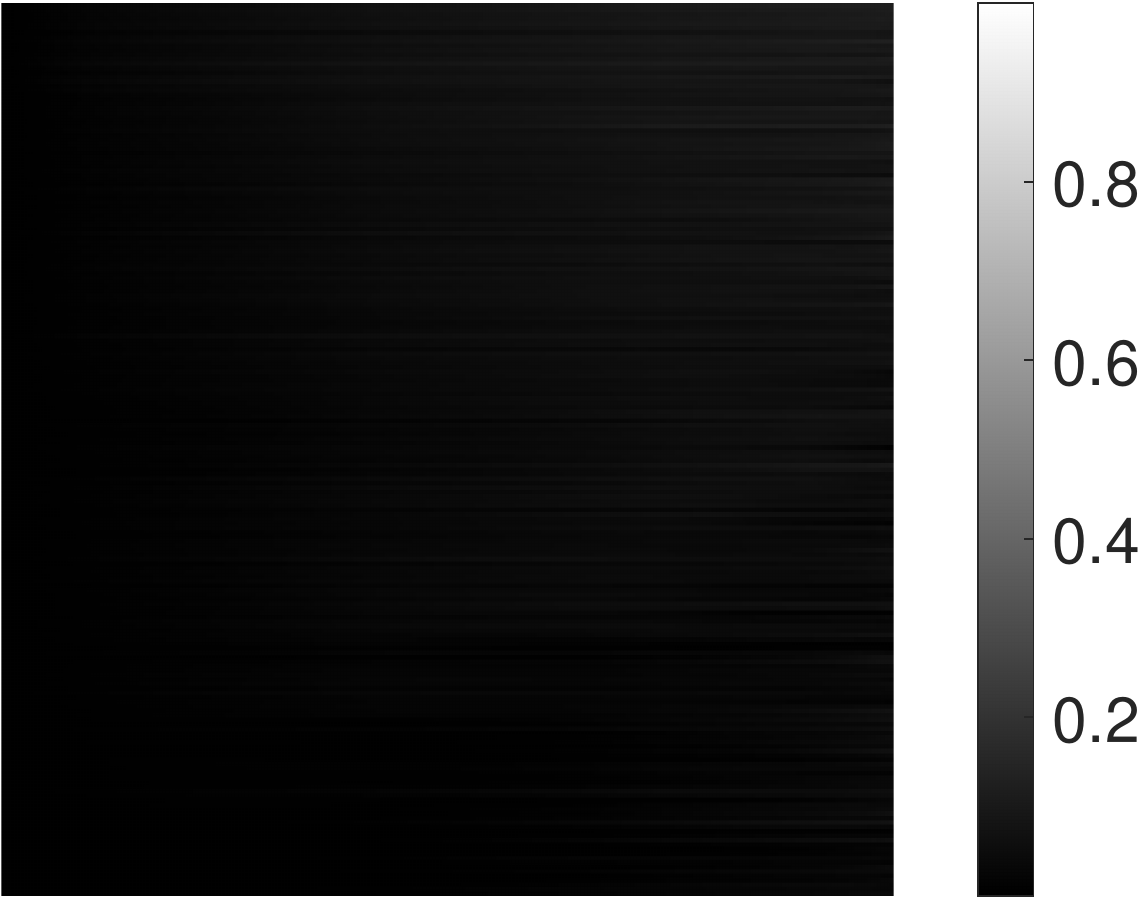} }
	\caption{Vector $\alpha_i$ (left panels), base-10 logarithmic plot of the corresponding scaled variances (middle panels) and vectors $\mW_i\alpha_i$ contributing to the final restoration (right panels) for $i=1$, i.e. vertical increments representation, (top row), $i=2$, i.e. horizontal increments representation, (bottom row).} 
	\label{fig:rec 2d}
\end{figure}

The restored image is shown in Figure \ref{fig:Data_2d}, while the contribution of the vertical and horizontal increment bases together with the output scaled variances corresponding to vectors $\alpha_1$ and $\alpha_2$ are shown in Figure \ref{fig:rec 2d}. We observe that the image is almost completely restored in terms of the basis vectors corresponding to increments in the vertical direction ($\alpha_1$), whereas the entries of $\alpha_2$, corresponding to increments in the horizontal direction is negligible. The representation in terms of $\mW_1$ is indeed sparser than that in terms of $\mW_2$, due to the shorter horizontal boundary of the white inclusion compared to the vertical boundary.

\medskip

\begin{figure}[!t]
	\centering
	\centerline{\includegraphics[width = 3.8cm]{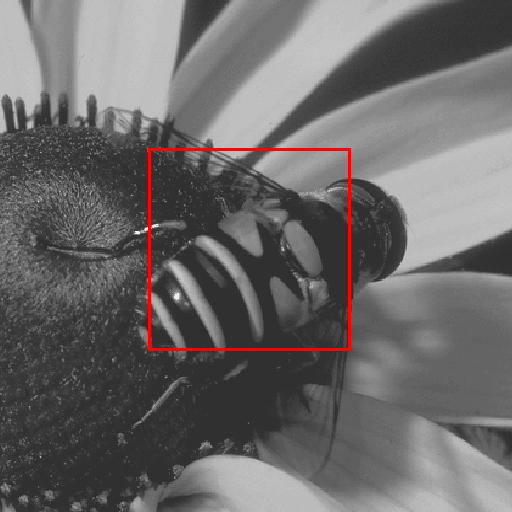}\quad \includegraphics[width = 3.8cm]{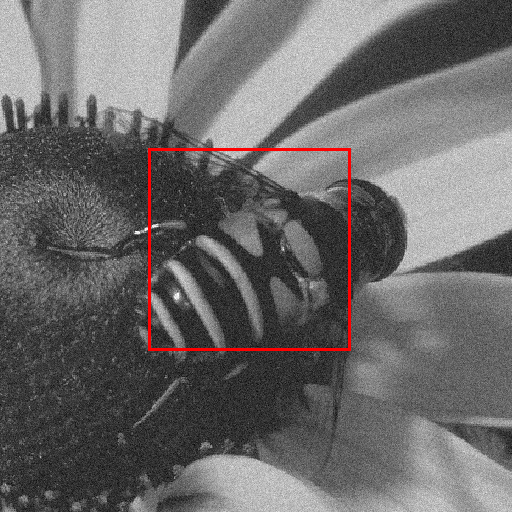}\quad\includegraphics[width = 3.8cm]{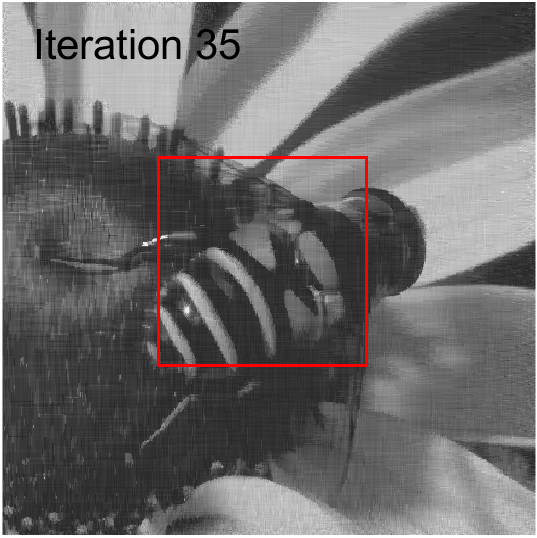}}
	\centerline{\includegraphics[width = 3.8cm]{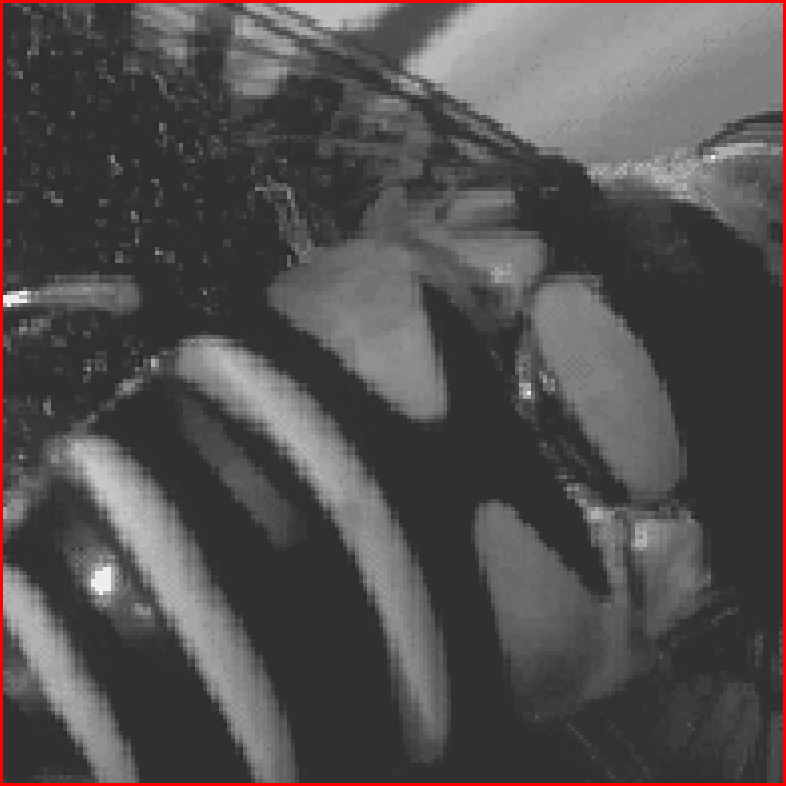}\quad\includegraphics[width = 3.8cm]{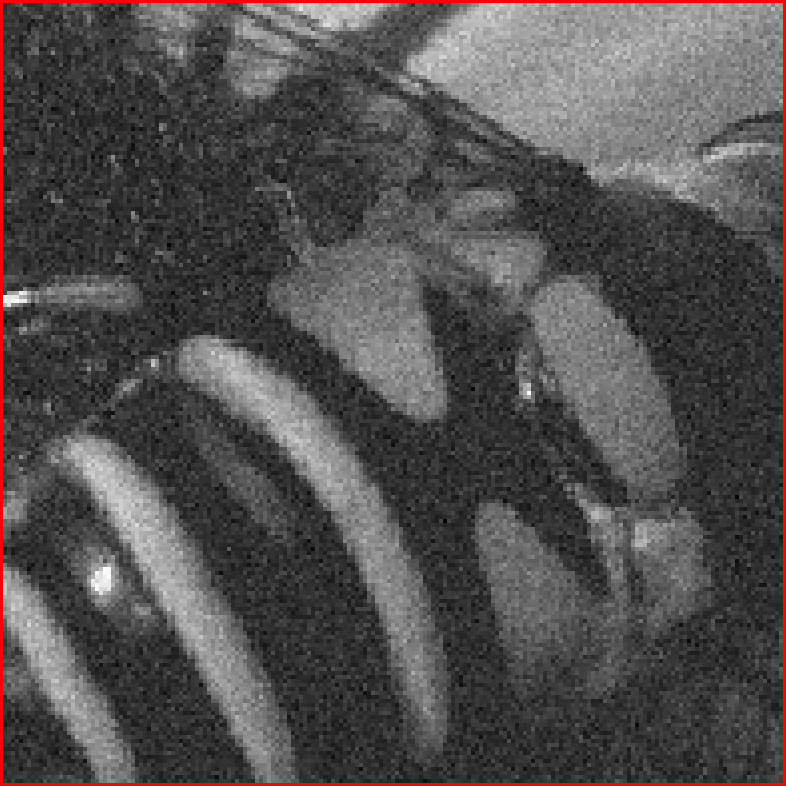}\quad\includegraphics[width = 3.8cm]{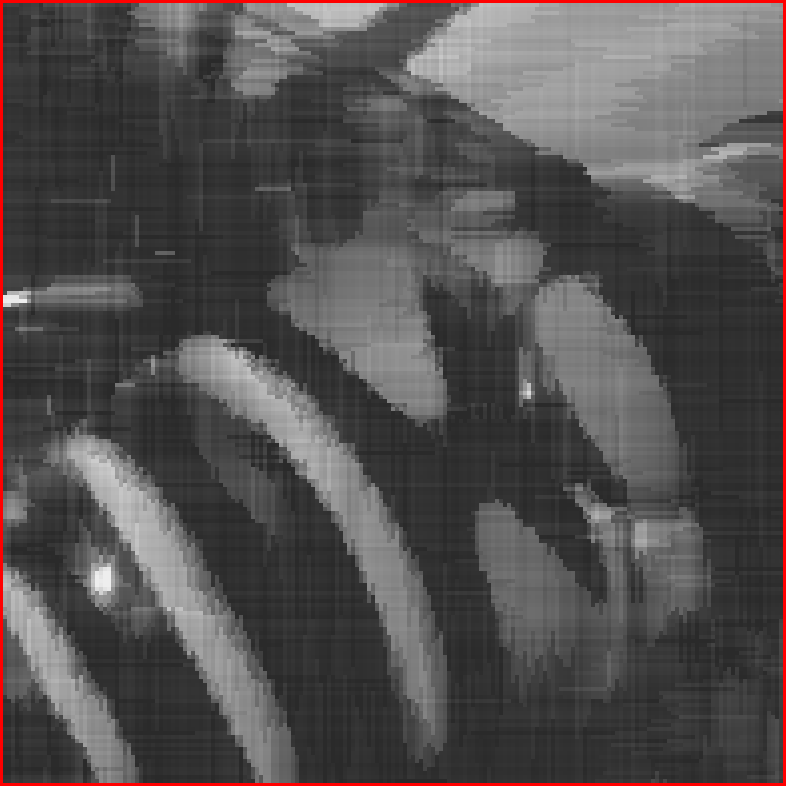}}
	\caption{\label{fig:bee_data}\emph{Top row}: original $512\times512$ test image (left), observed data (middle) and denoised image (right). \emph{Bottom row}: respective close-up(s).}
\end{figure}

	\paragraph{\textbf{Image denoising on a natural image}}
	We demonstrate the scalability of our approach with large-scale denoising problem. Consider the $512 \times 512$ gray-scale natural image with pixel values between $0$ and $1$ shown in the top left panel of Figure \ref{fig:bee_data}. The observed data $b$ is a version corrupter by added white Gaussian noise with standard deviation $\sigma$ set to $5\%$ of the maximum of the noiseless image, i.e. $\sigma=0.05$,  shown in the top middle panel of Figure \ref{fig:bee_data}.

	We consider an overcomplete basis $\mathsf{W}=[\mathsf{W}_1,\,\mathsf{W}_2]\in\R^{n^2\times 2n^2}$, with $\mathsf{W}_1,\mathsf{W}_2\in\R^{n^2\times n^2}$, $n=512$, defined in \eqref{eq:dict 2d}; the vector describing the 2D signal in the selected dictionary can thus be written as $\alpha=[\alpha_1,\,\alpha_2]^\mT$, with $\alpha_1,\alpha_2\in\R^{n^2}$ representing the vertical and horizontal increments, respectively.
	Note that since this test image is not piecewise constant, but rather a mixture of jumps, smooth, and textured parts, we do not expect the vertical and horizontal increments to be \emph{naturally} sparse. In fact, when representing the original image in the chosen dictionary, about 47\% of the coefficients $\alpha_1,\alpha_2$ in both horizontal and vertical directions are non-zero.

We run the global hybrid IAS with $(r^{(1)},\eta^{(1)})=(1,10^{-4})$, $(r^{(2)},\eta^{(2)})=(1/2,10^{-3})$, and letting  $\vartheta^{(1)},\vartheta^{(2)}$ be set automatically.

\begin{figure}
	\centerline{\includegraphics[height =3.7cm]{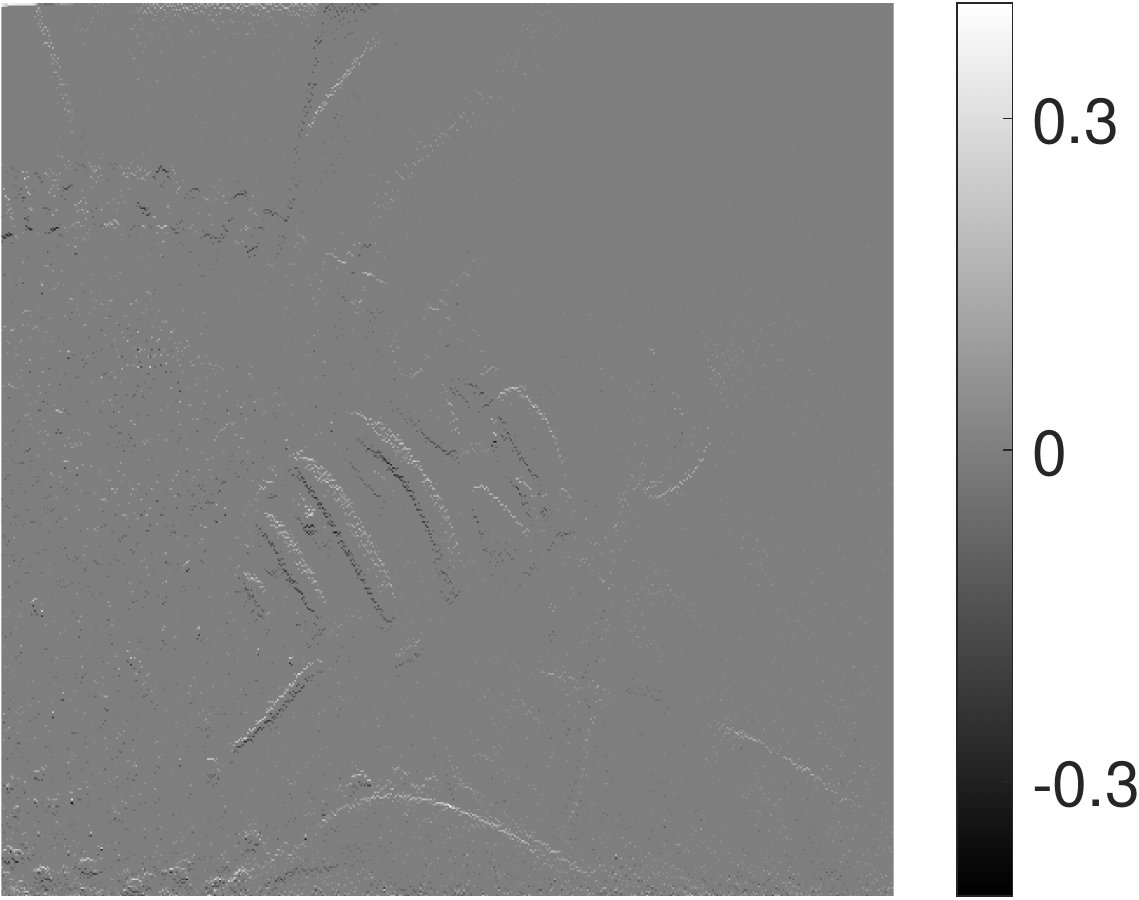}\;\includegraphics[height = 3.7cm]{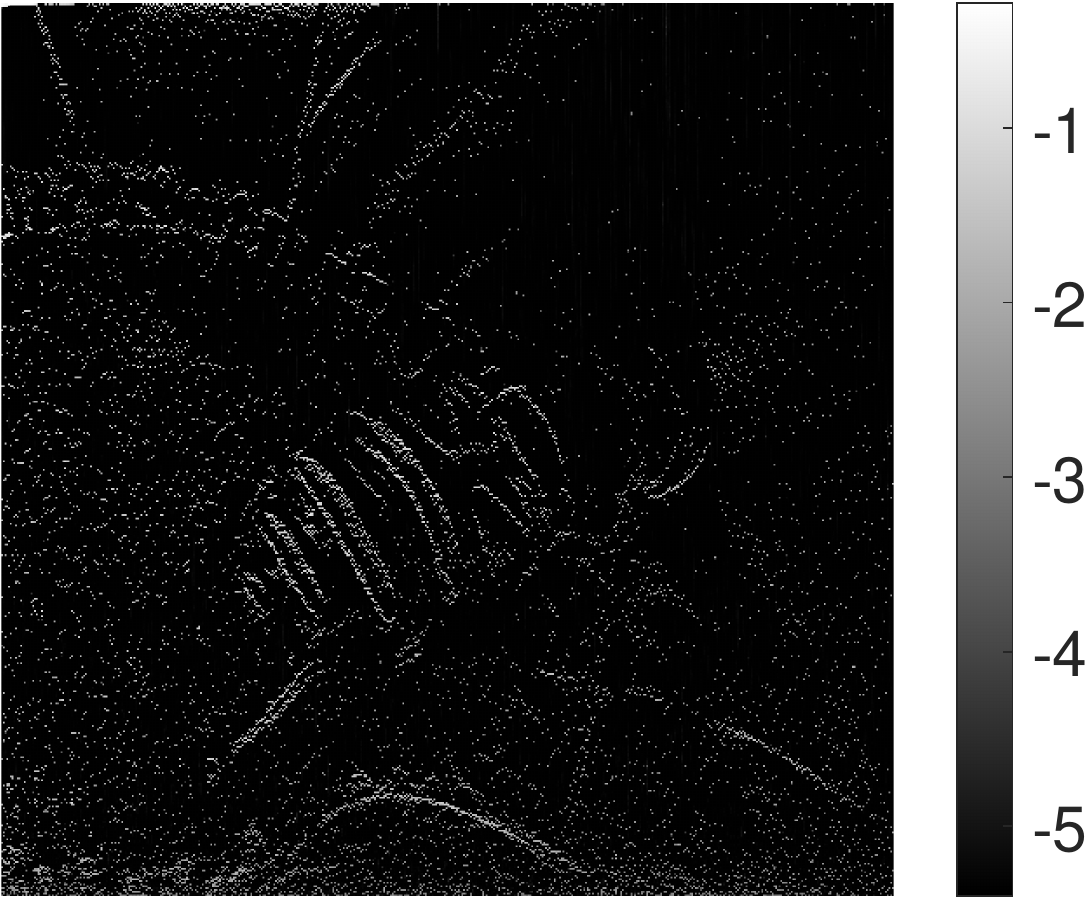}\;\includegraphics[height = 3.7cm]{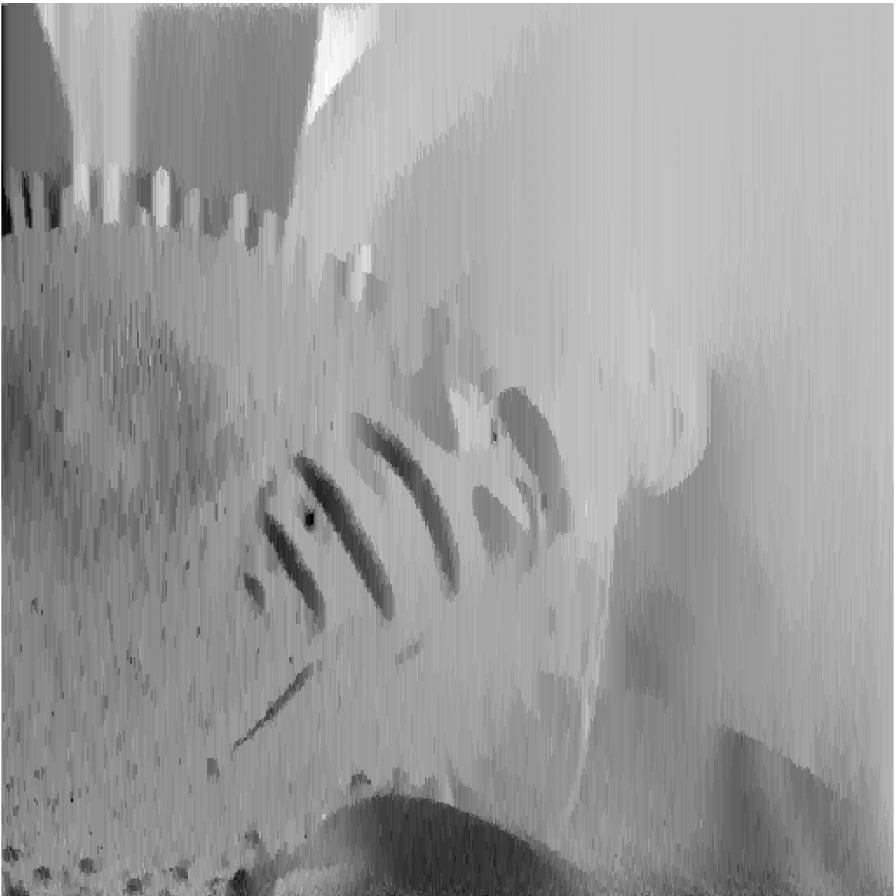}}
	\centerline{\includegraphics[height =3.7cm]{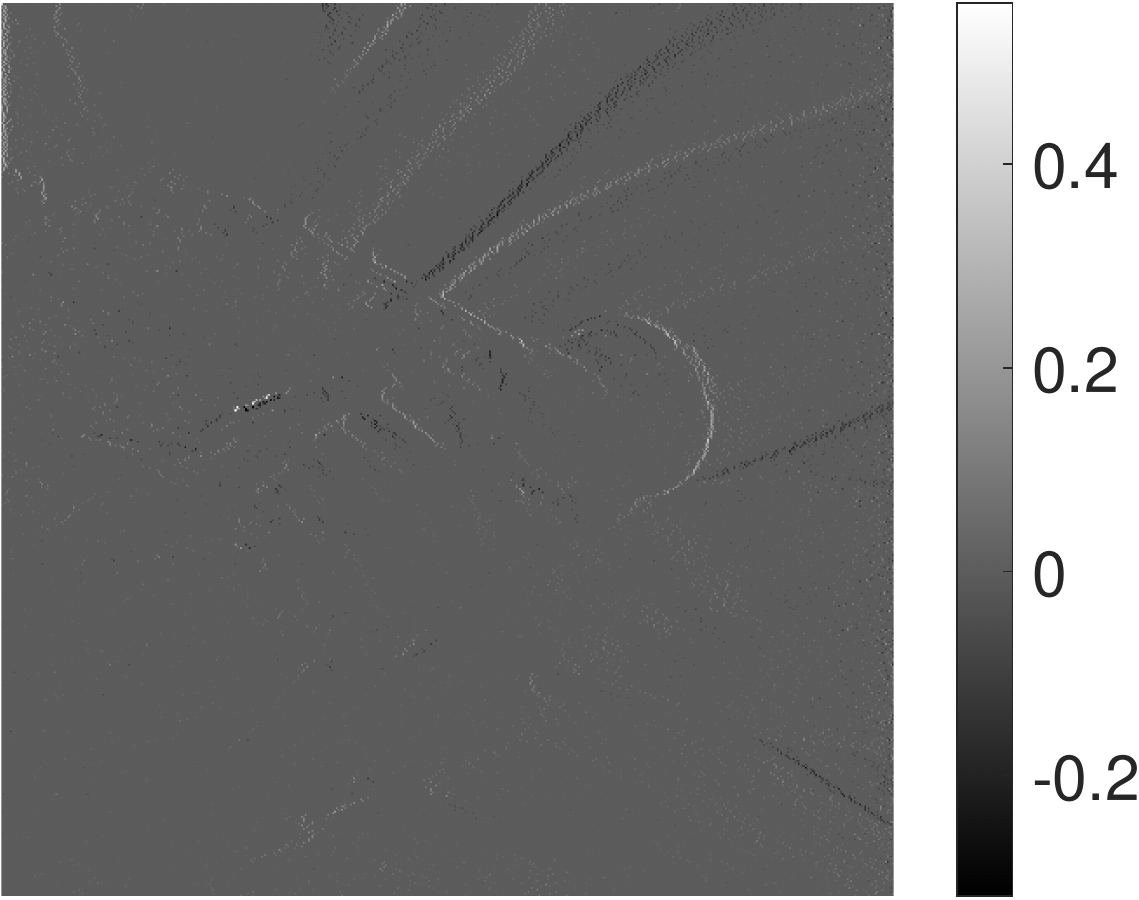}\;\includegraphics[height = 3.7cm]{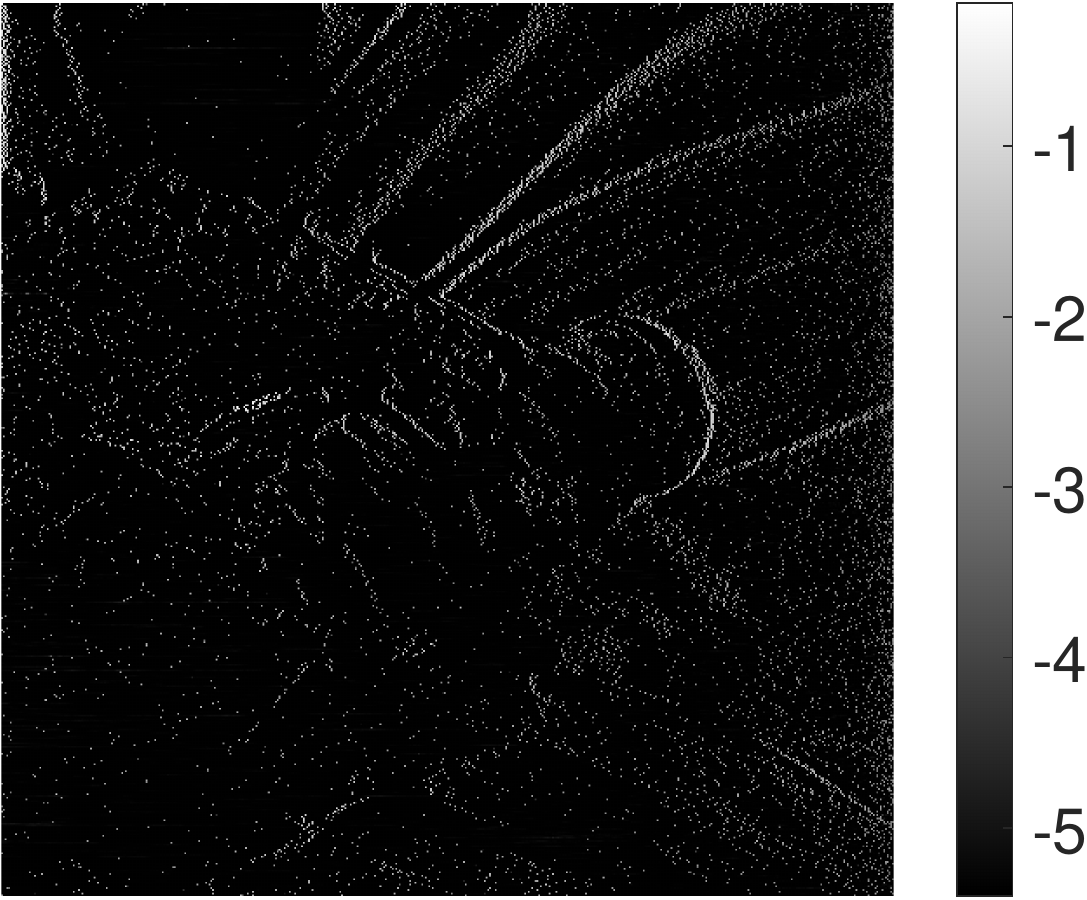}\;\includegraphics[height = 3.7cm]{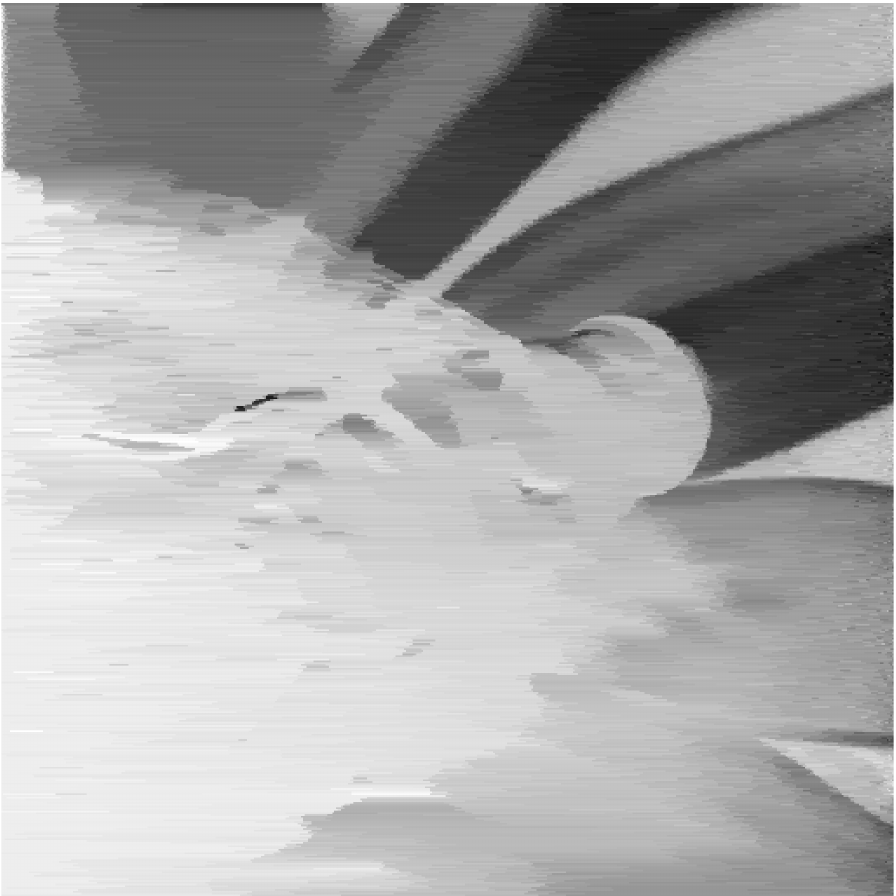}}
	\caption{Vector $\alpha_i$ (left panels), base-10 logarithmic plot of the corresponding scaled variances (middle panels), and vectors $\mathsf{W}_i\alpha_i$ (right panels) for $i=1$, i.e. vertical increments representation, (top row) and $i=2$, horizontal increments representation, (bottom row).} 
	\label{fig:bee_th}
\end{figure}

The restored image is shown in the top right panel of Figure \ref{fig:bee_data}, together with a close-up to facilitate the visual comparison. In Figure \ref{fig:bee_th}, we also show the output representation vectors $\alpha_1,\alpha_2$, with the corresponding scaled output variances and their contribution in the final restoration.

\begin{figure}
	\centerline{\includegraphics[height =5.2cm]{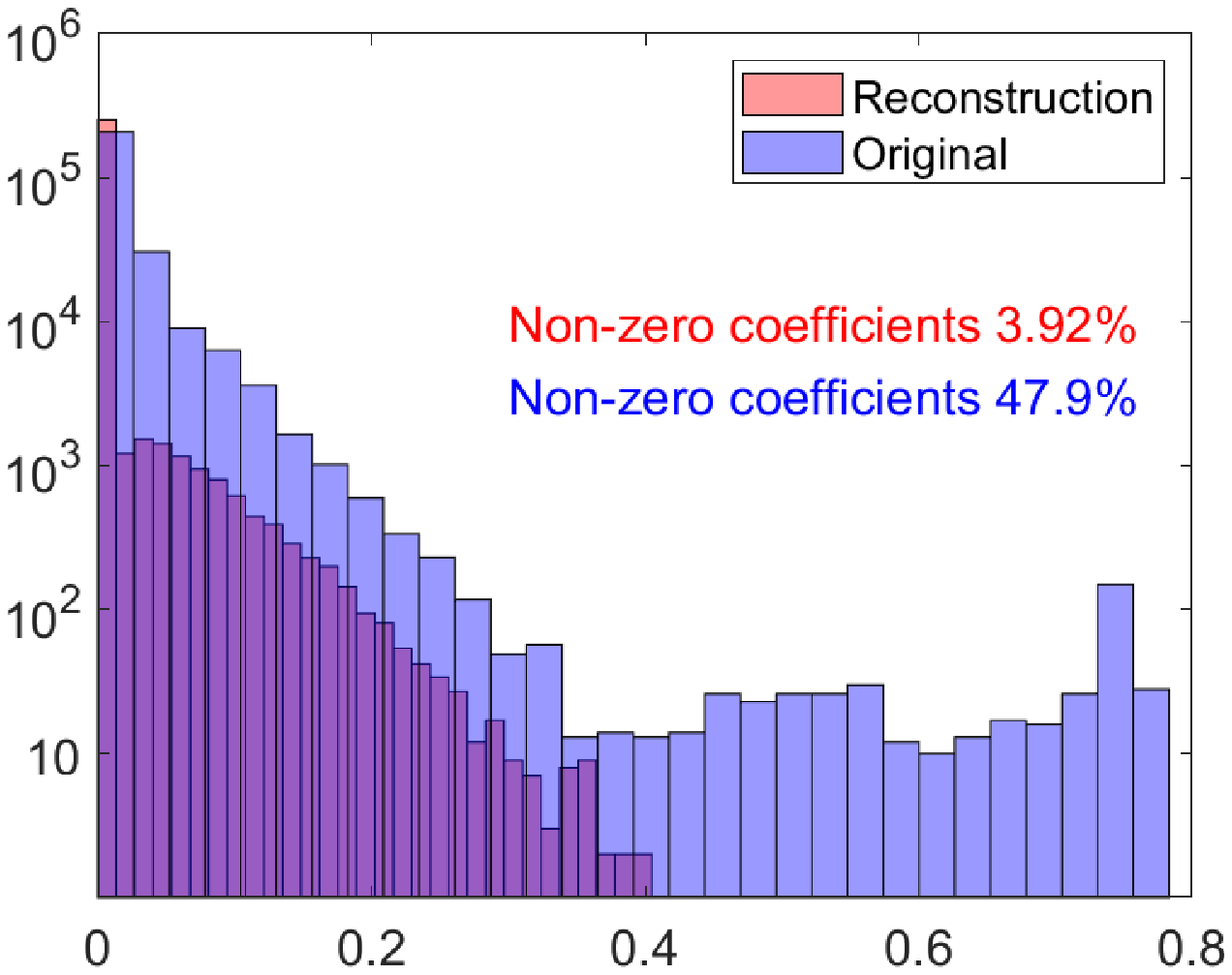}\;\includegraphics[height = 5.2cm]{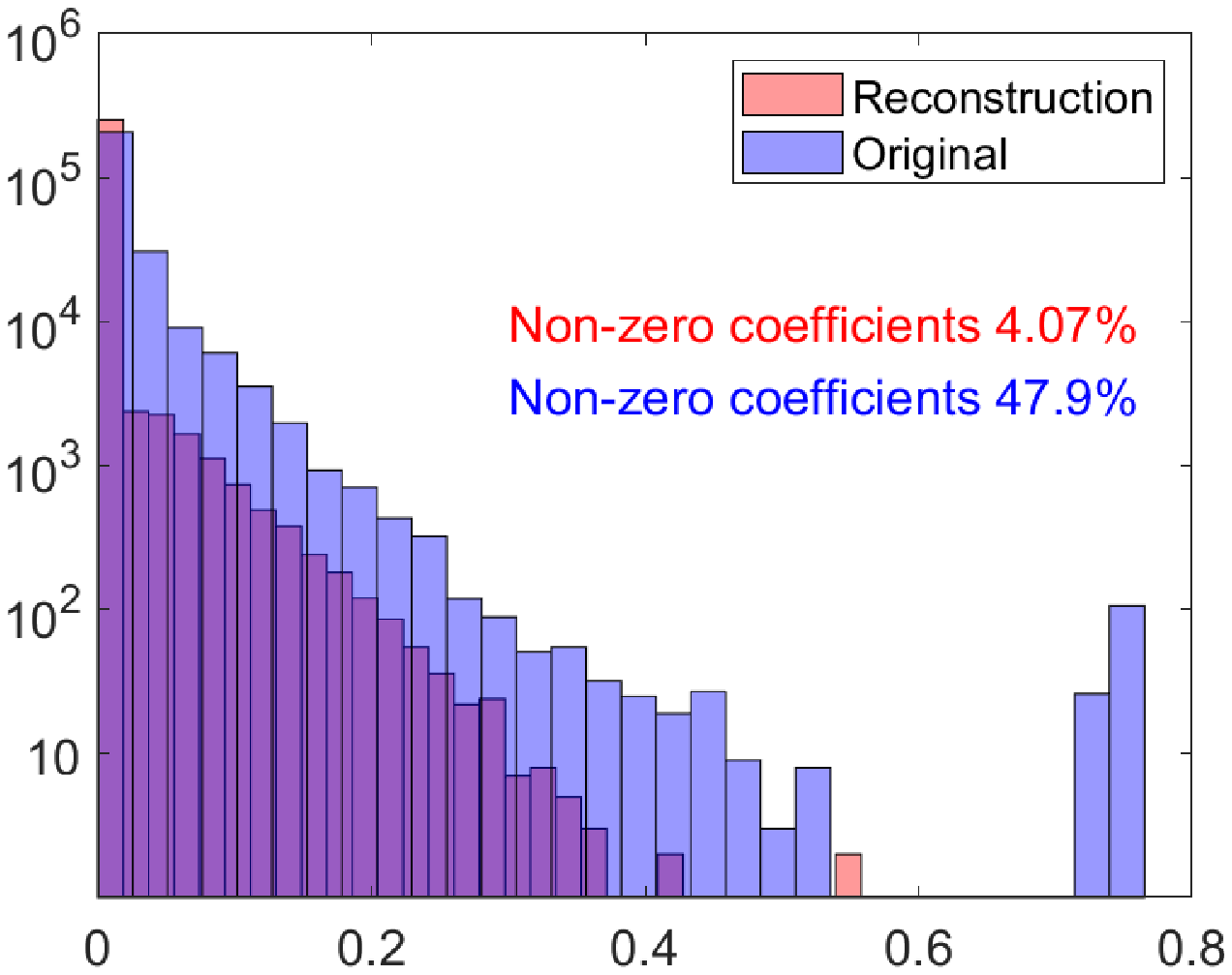}}
	\caption{Value distributions of the coefficients corresponding to vertical (left) and horizontal (right) increments in logarithmic scale. The blue distributions correspond to the original image which does not allow sparse representation in the basis, which is reflected in the high percentage of non-vanishing coefficients, while the red distribution is the sparse denoised reconstruction, in which the percentage of coefficients above the negligible threshold value is reduced by an order of magnitude.}
	\label{fig:bee_hist}
\end{figure}

	Besides the quality of the restored image, we are interested in highlighting the compressing capability of our approach. Consider first the coefficients $\alpha=[\alpha_1,\,\alpha_2]^\mT$ of the original image, and choose a threshold value $\beta>0$ to be the smallest non-zero coefficient,
	\begin{equation}\notag
	\beta = \min\{|\alpha_j| \mid |\alpha_j|>0\}.
	\end{equation}
	After computing the restored coefficients using the sparsity promoting hybrid IAS, setting to zero those whose absolute value is below a threshold value, we find that only 3.92\% of those for the vertical components, and 4.07\% of those for the horizontal components are nonvanishing. Thus, the representation of the image in this basis was compressed by a factor more than ten from the original image without a significant deterioration in the image quality. Figure~\ref{fig:bee_hist} shows the value distributions of the original and restored coefficients, plotted as histograms in logarithmic scale: clearly the coefficient values are significantly compressed towards zero in both directions.  It is worth remarking here that the Bayesian target in sparsity promoting problems has to be understood as compressibility, as the entries of $\alpha$ cannot vanish by construction of the algorithm, but they can be made arbitrary small by suitable parameter selection.



\begin{figure}[!t]
	\centering
	\centerline{
		\includegraphics[height=3.75cm]{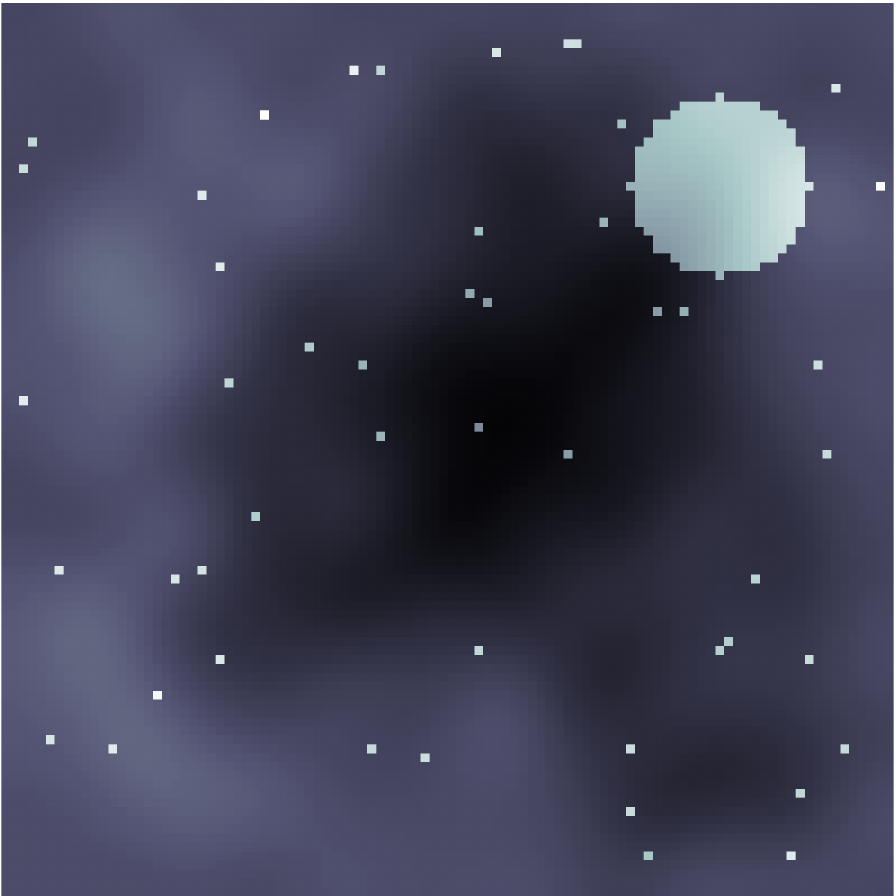}\;
		\includegraphics[height=3.75cm]{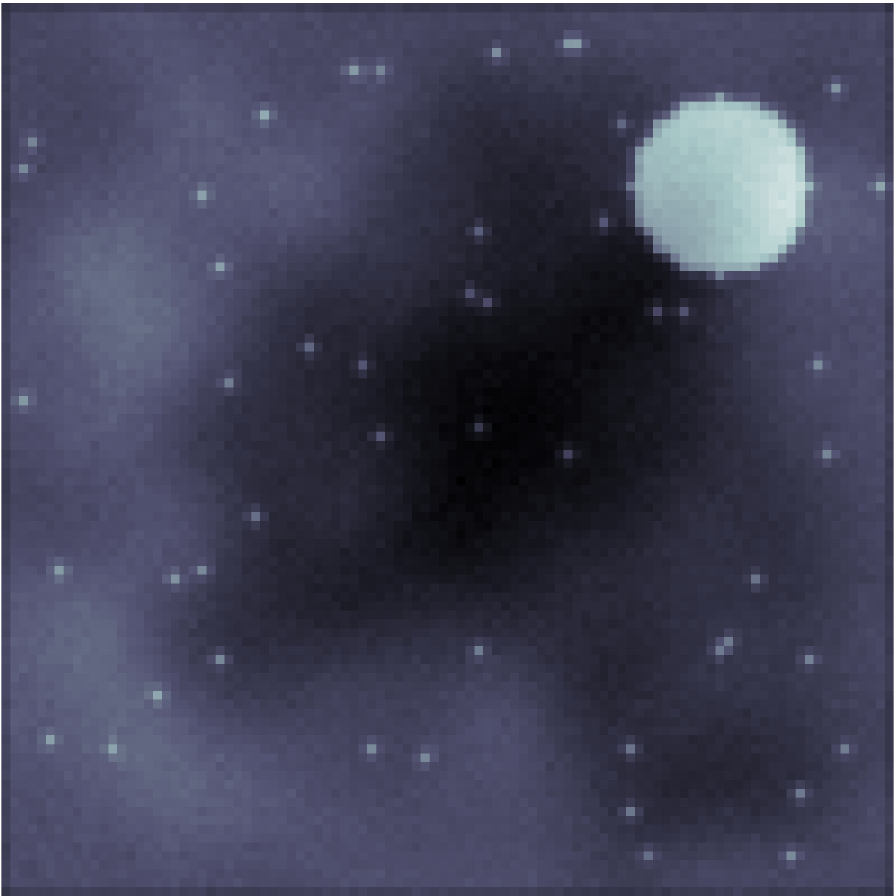}\;	\includegraphics[height=3.75cm]{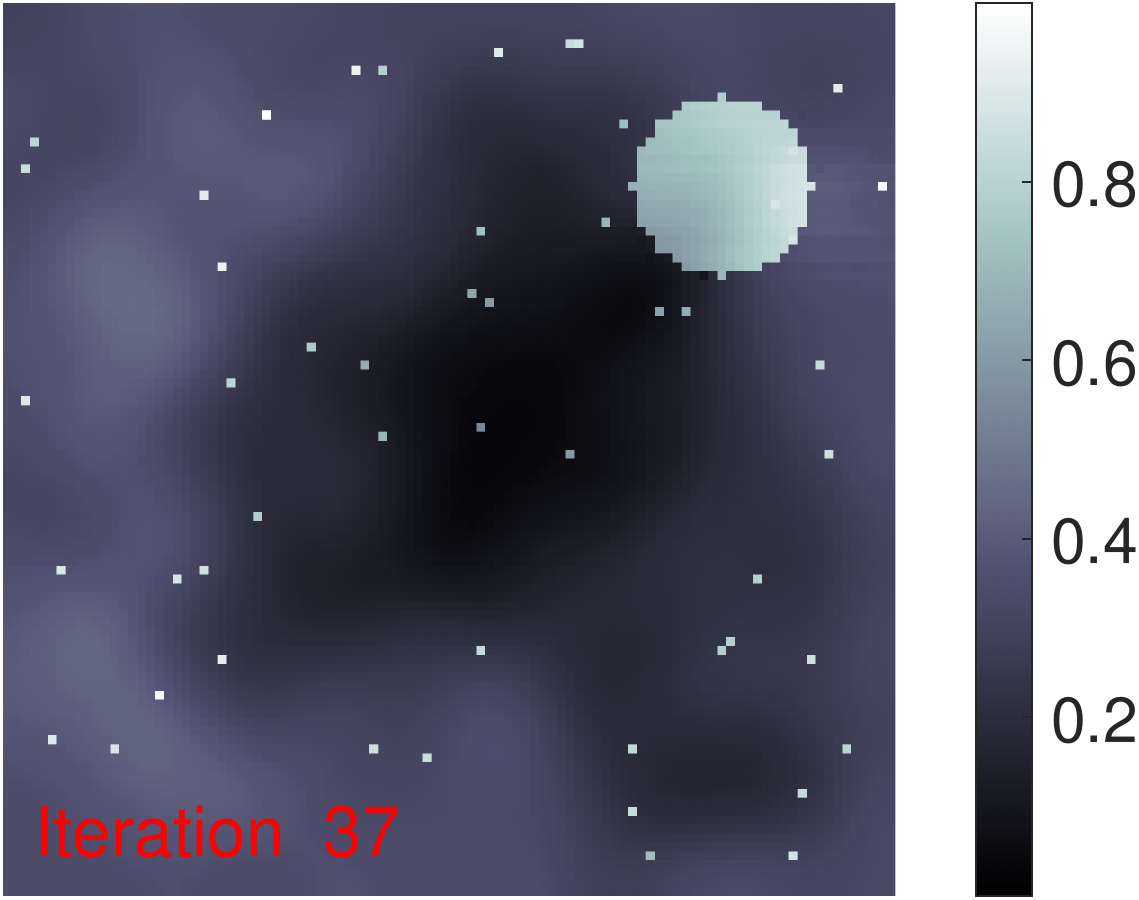}	}
	\caption{Original image (left), observed data (middle) and reconstructed image (right).}
	\label{fig:Data_2d_2}
\end{figure}

\paragraph{\textbf{Image restoration}} In the fourth example, we consider the restoration problem of the $n \times n$ generative image in Figure \ref{fig:Data_2d_2}, with $n=100$,  with values in $[0,1]$. The image has been corrupted by Gaussian blur of width $w=0.006$ and additive scaled white Gaussian noise with standard deviation $\sigma$ set to $1\%$ of the maximum of the noiseless signal, i.e. $\sigma=0.01$ - see Figure \ref{fig:Data_2d_2}. The test image presents three distinctive features, namely point-wise stars, the blocky moon and the smooth cloud. After re-arranging the original $x$ in a vectorized form by stacking its entries in columnwise order, we hypothesize that a suitable dictionary for the problem of interest is $\mW = [\mW_1 \,,\, \mW_2 \,,\, \mW_3 \,,\, \mW_4] \in \R^{n^2 \times 4n^2}$, with
\begin{equation}
\mW_1 {=} \mI_{n^2}\,,\;\mW_2 {=} (\mI_n \otimes \mB)^{-1} \in \R^{n^2 \times n^2}\,,\; \mW_3 {=} (\mB \otimes \mI_n)^{-1} \in \R^{n^2 \times n^2} \text{ and } \mW_4 {=} \mC^\mT\,,
\end{equation}
where $\mB$ is as in \eqref{difference} and $\mC$ is the 2D cosine transform matrix. The problem is to estimate the sparse vector $\alpha = [\alpha_1\,,\,\alpha_2\,,\,\alpha_3\,,\,\alpha_4]^{\mT}$, with $\alpha_i\in\R^{n^2}$, for $i=1,2,3,4$, from the data vector $b\in\R^{n^2}$, given the forward model
\begin{equation}
b = \mA\mW \alpha + \varepsilon =  [\mW_1 \,,\, \mW_2 \,,\, \mW_3 \,,\, \mW_4] \begin{bmatrix}
\alpha_1\\
\alpha_2\\
\alpha_3\\
\alpha_4
\end{bmatrix}+\varepsilon, \quad \varepsilon \sim \mathcal{N}(0,\sigma^2 \mI_{n^2})\,,
\end{equation}
with $\mA$ representing the discrete blur operator.

The global hybrid IAS is run with hyperparameters $(r^{(1)},\eta^{(1)})=(1,10^{-4})$, $(r^{(2)},\eta^{(2)})=(1/2,10^{-4})$ and $\vartheta^{(1)},\vartheta^{(2)}$ automatically fixed as in the previous examples.

\begin{figure}[!t]
	\centering
	\centerline{\includegraphics[height=3.75cm]{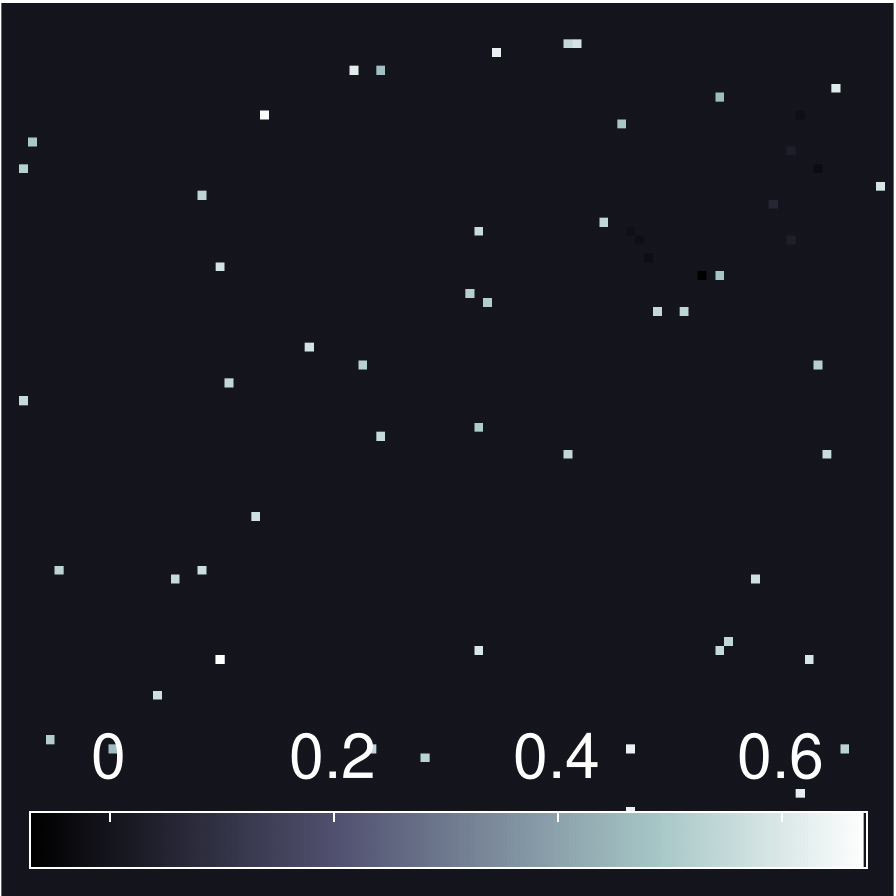}\;	\includegraphics[height=3.75cm]{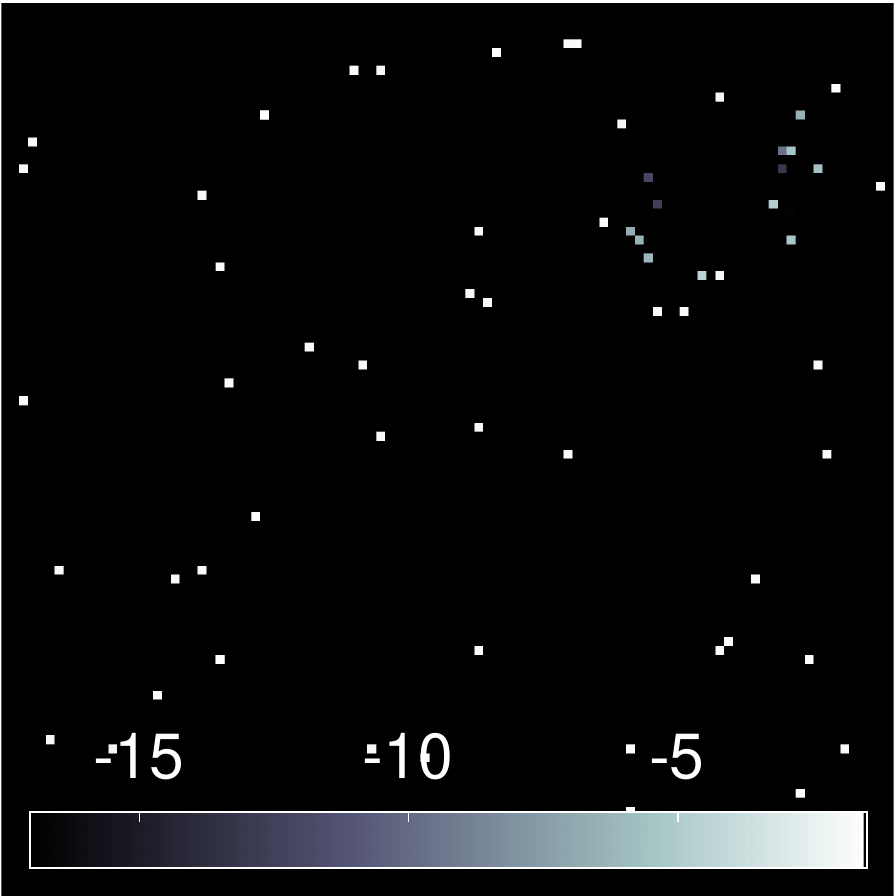}\; 		\includegraphics[height=3.75cm]{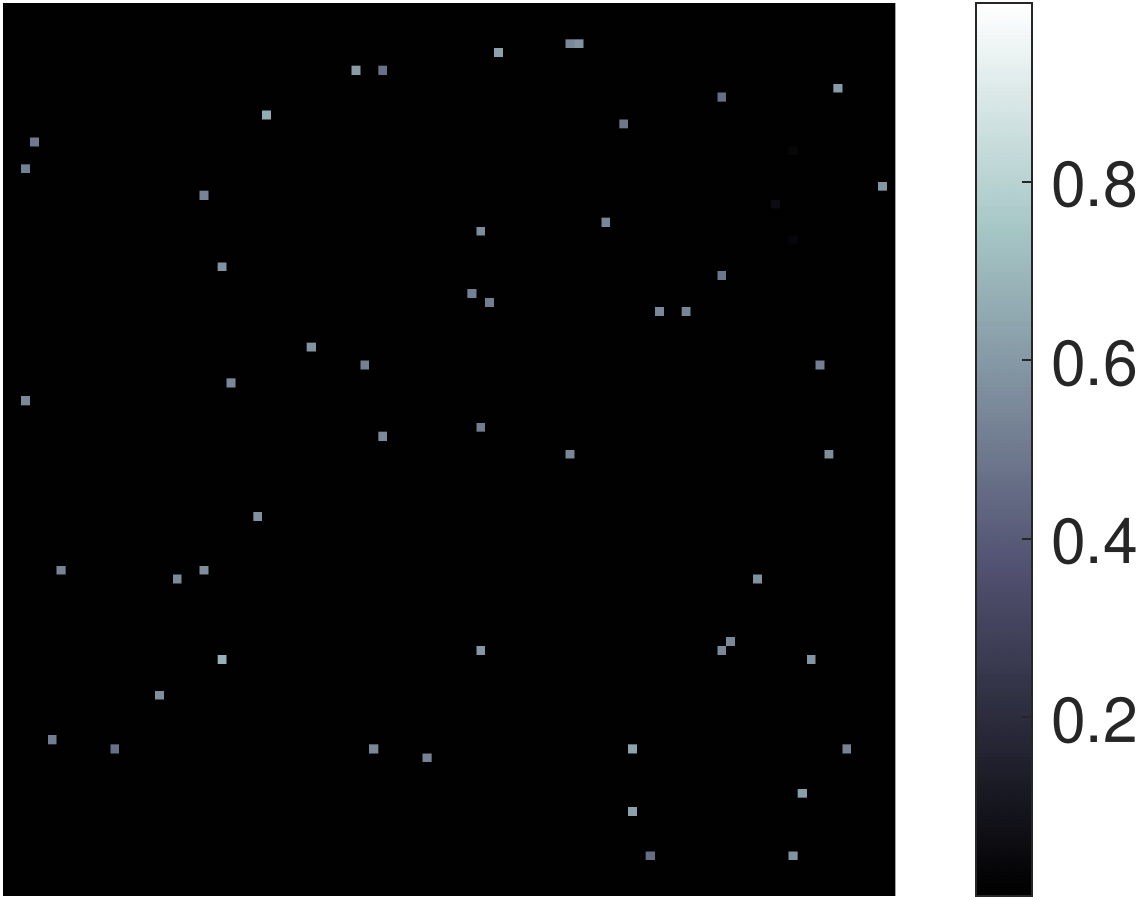} }
	\centerline{\includegraphics[height=3.75cm]{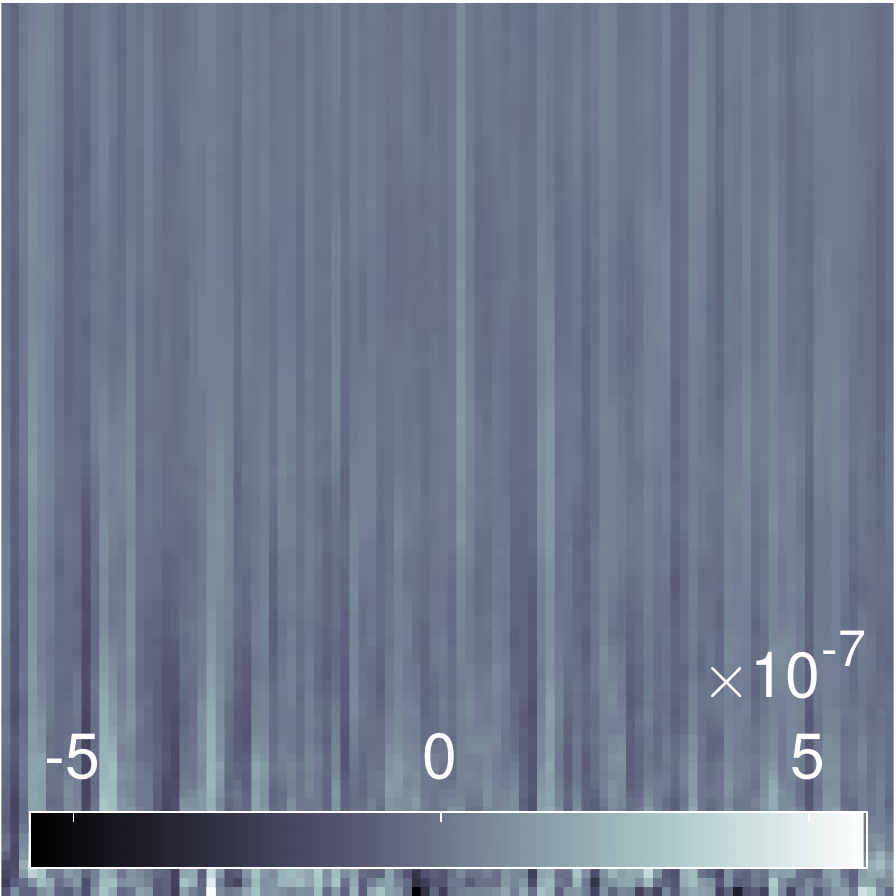}\;		\includegraphics[height=3.75cm]{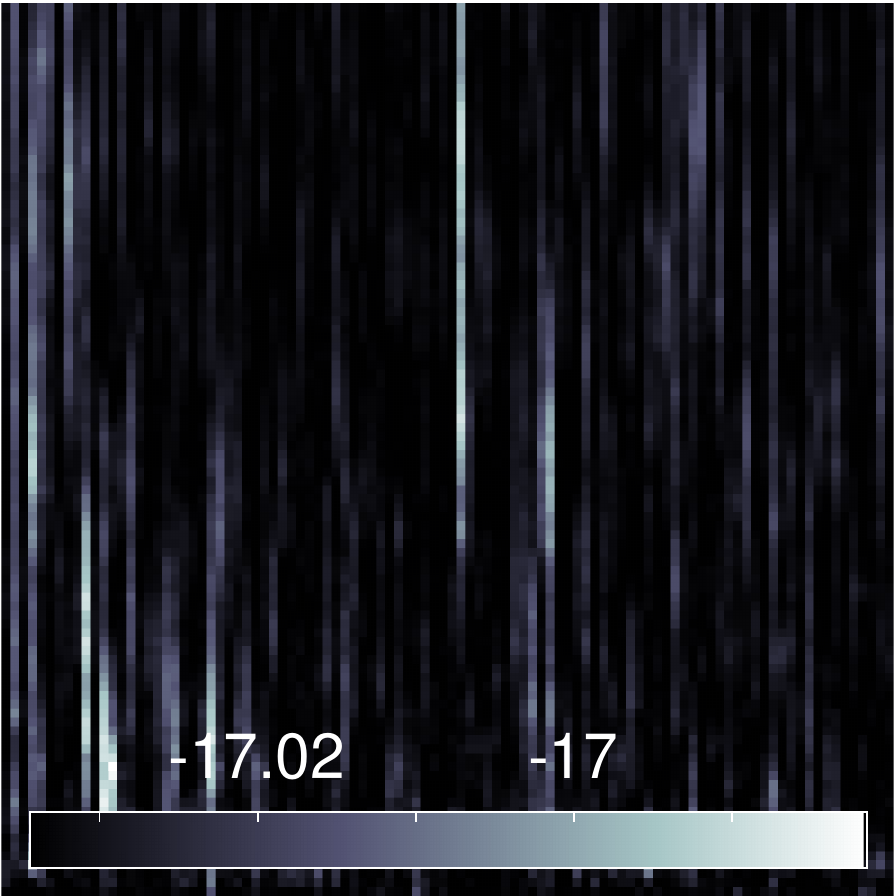}\; 		\includegraphics[height=3.75cm]{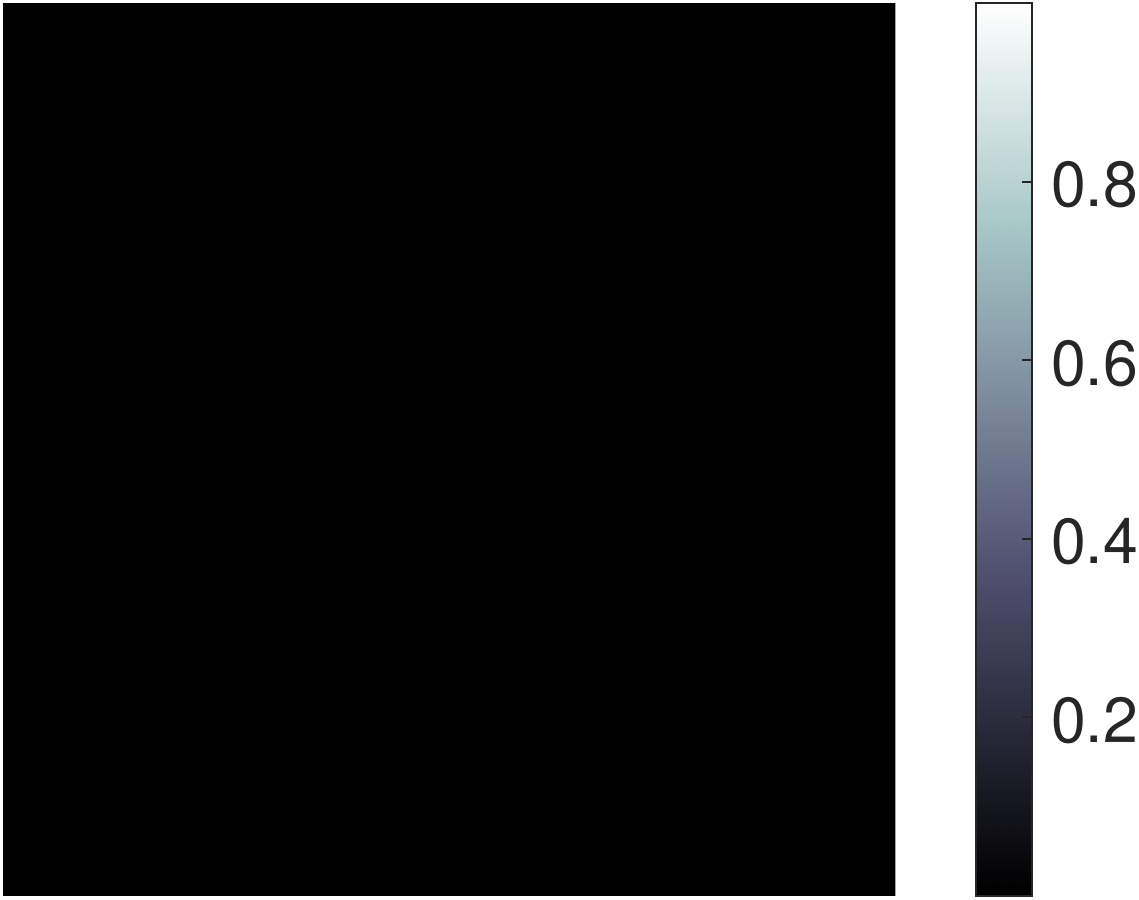} }	
	\centerline{\includegraphics[height=3.75cm]{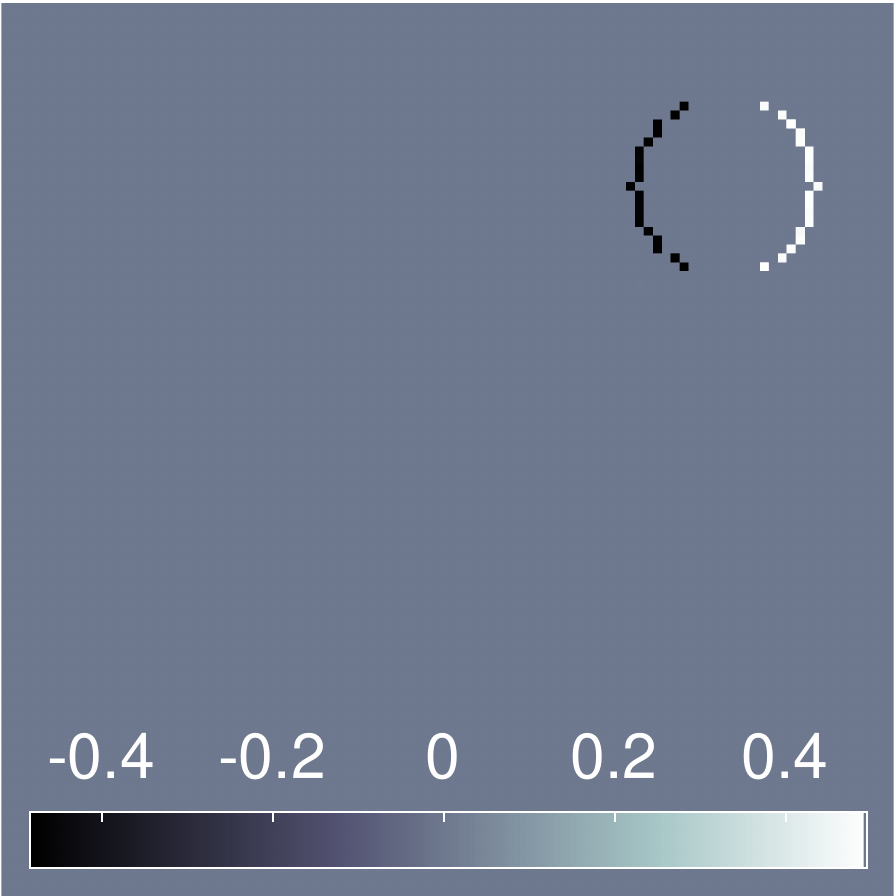}\;		\includegraphics[height=3.75cm]{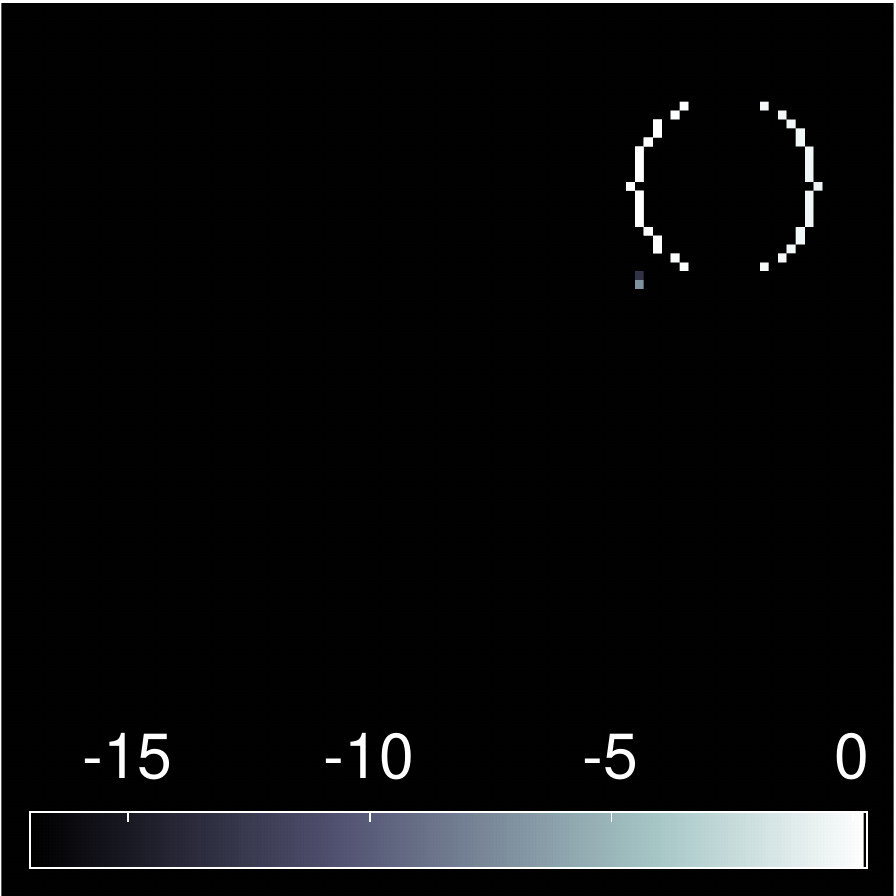}\; 		\includegraphics[height=3.75cm]{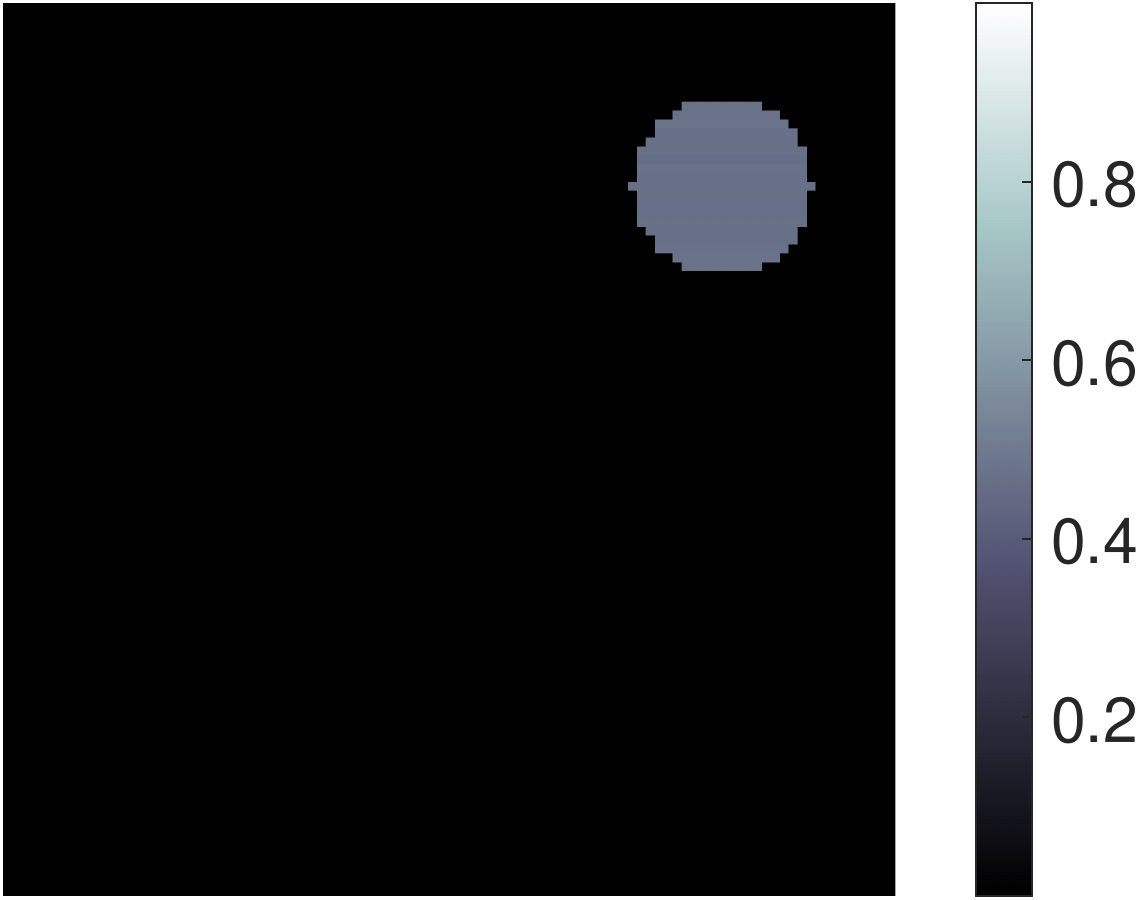} }	
	\centerline{\includegraphics[height=3.75cm]{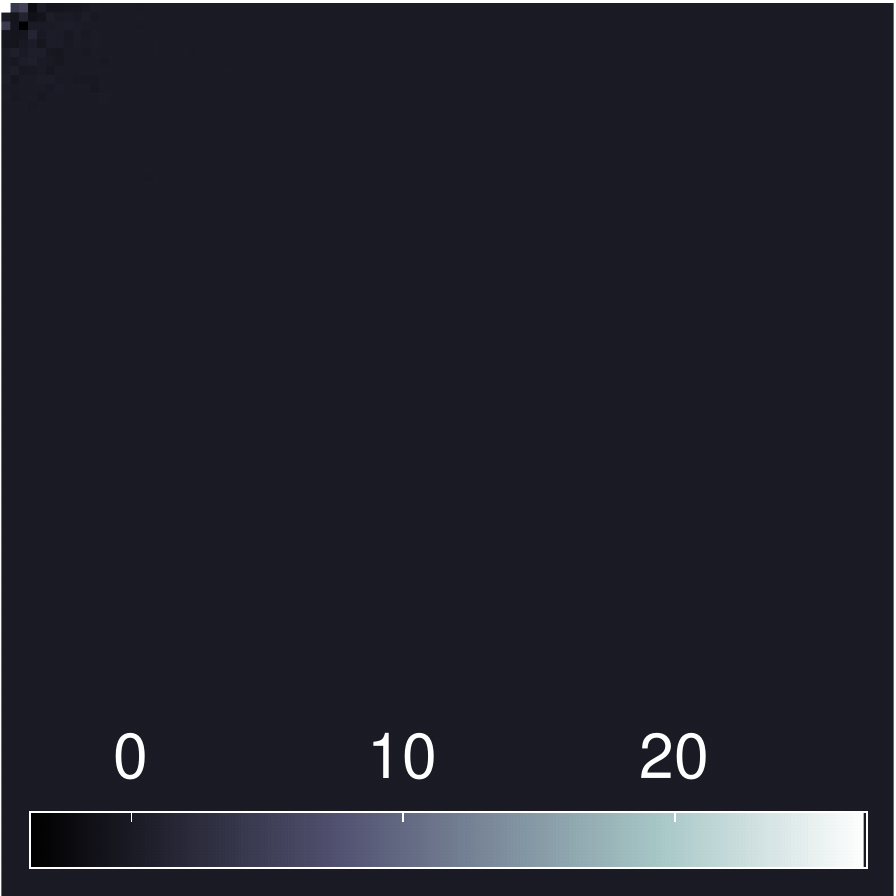}\;		\includegraphics[height=3.75cm]{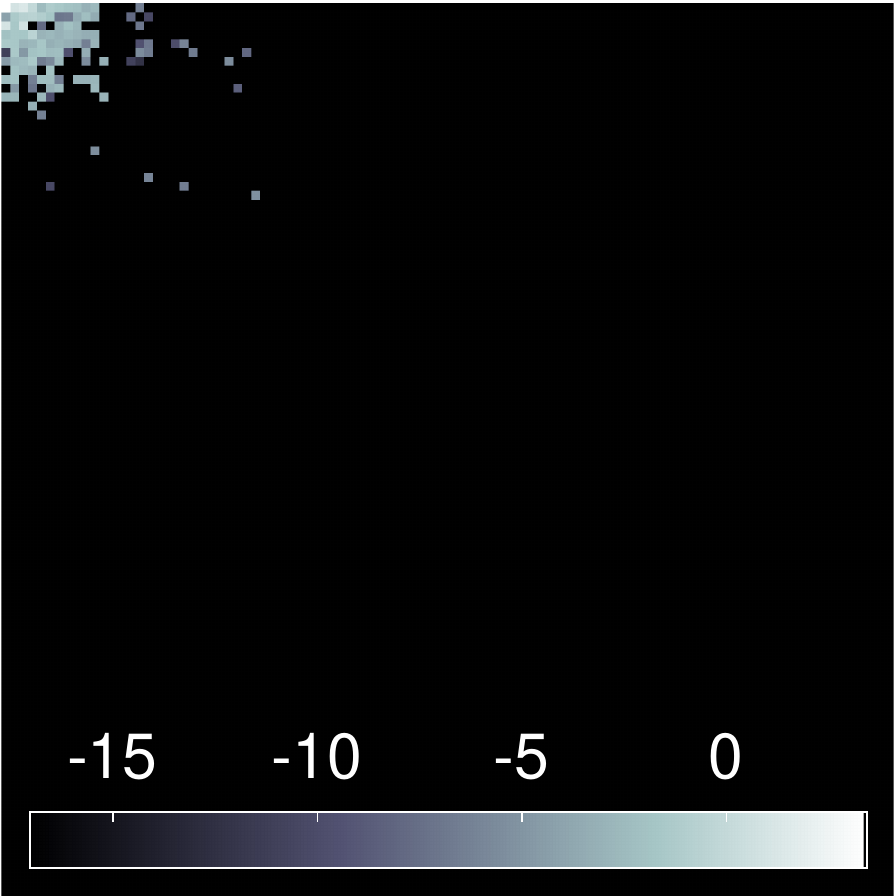}\; 		\includegraphics[height=3.75cm]{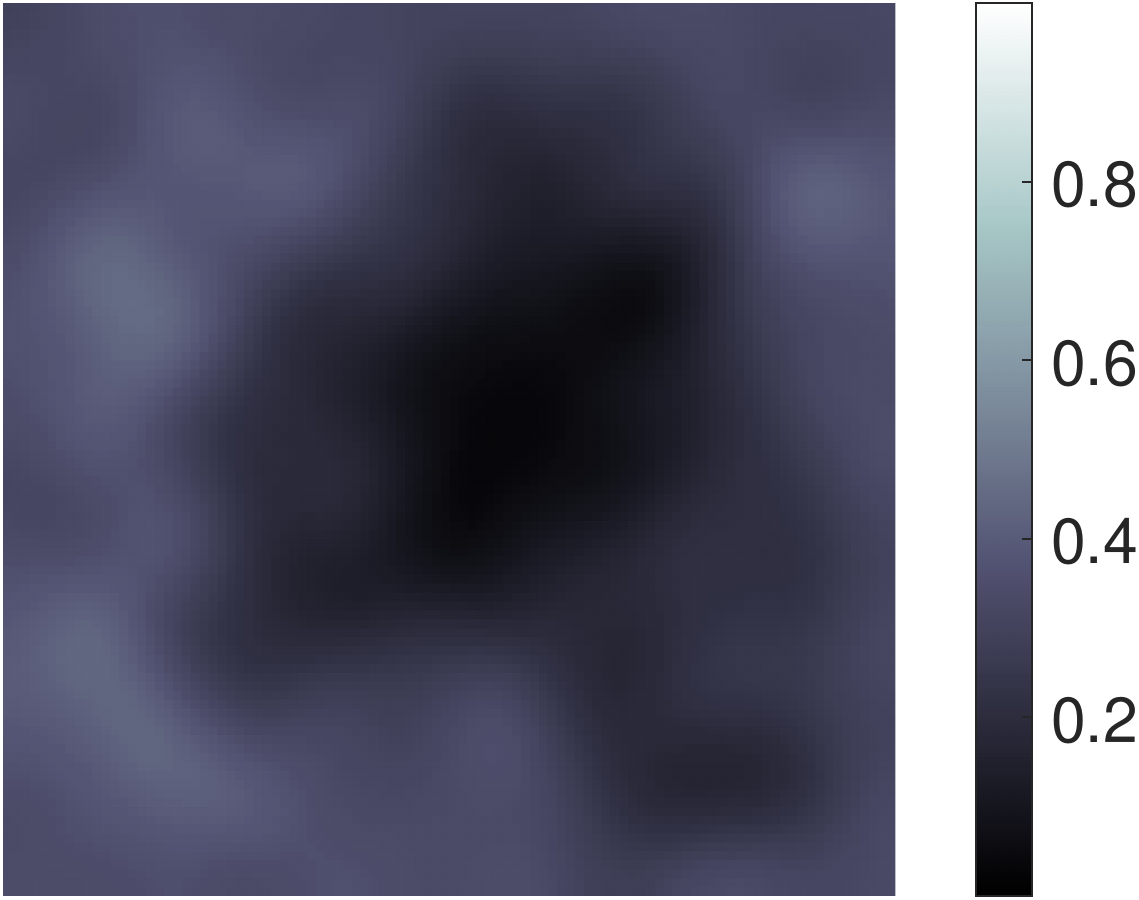} }
	\caption{Representation vectors $\alpha_i$ (left panels), base-10 logarithmic plot of the corresponding scaled variances (middle panels) and vectors $\mW_i\alpha_i$ contributing to the final restoration (right panels) for $i=1$ (first row), $i=2$ (second row), $i=3$ (third row) and $i=4$ (fourth row). }
	\label{fig:rec 2d 2}
\end{figure}

The image restored via the global hybrid IAS algorithm is shown in Figure \ref{fig:Data_2d_2}, while Figure \ref{fig:rec 2d 2} shows the reconstructions of the representation vectors $\alpha_i$, the corresponding variances scaled by the sensitivities, and the contribution of the vectors $\mW_i\alpha_i$ in the final restoration, for $i=1,2,3,4.$

We point out that, as in the previous example, the representation vectors in both the vertical and horizontal increment bases are sparse. Nonetheless, the hybrid hypermodel selects the one with fewer non-zero entries.



\begin{figure}
	\centerline{\includegraphics[width=10cm]{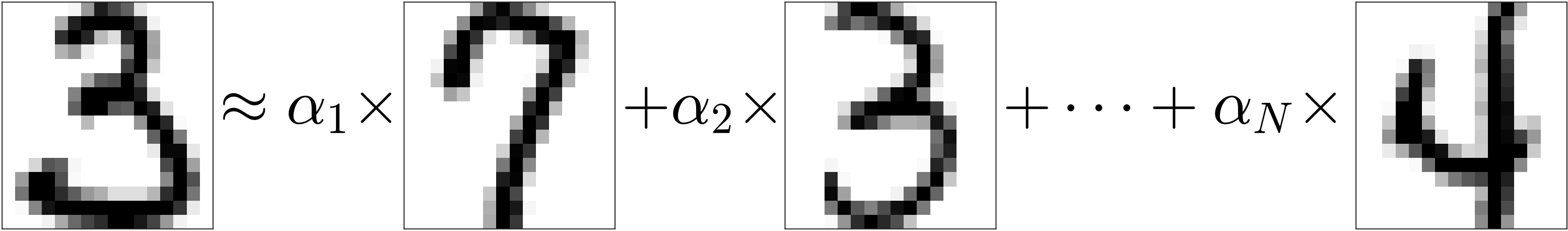}}
	\caption{\label{fig:DLcartoon} A schematic representation of the dictionary learning example. The digit on the right is the non-annotated image $b$, which is approximated in terms of the annotated atoms $w_i$ on the right. The coefficients $\alpha_i$ can then be used to identify the digit.}
\end{figure}

\medskip
\paragraph{\textbf{Dictionary learning}} The final example, coming from machine learning, is concerned with the sparse identification of hand-written digits based on a dictionary of annotated data. 
Consider the MNIST data set of hand-written digits $0,1,\ldots,9$ digitized as $16\times 16$ black-and-white images. Denoting by $w^{(j)}\in\R^{(16)^{2}}$, $1\leq j\leq N$ the vectorized image vectors of $N = 1\,707$ handwritten digits constituting the atoms of the dictionary, and by $c_j\in\{0,1,\ldots,9\}$ the corresponding annotations, we form the dictionary matrix
\begin{equation}
\mW = \left[\begin{array}{ccc} w^{(1)} &\cdots & w^{(N)}\end{array}\right] \in \R^{n\times N}, \quad n=256,\; N= 1\,705.
\end{equation}
To identify an handwritten digit $b$ drawn from an independent set of handwritten digits, we seek to represent it in a sparse manner in terms of the given dictionary,
\begin{equation}
b = \mW \alpha +\varepsilon,
\end{equation}
where $\alpha\in\R^N$ is a sparse vector, and $\varepsilon$ represents the discrepancy between the data and its representation. The idea is represented schematically in Figure~\ref{fig:DLcartoon}.
We point out that in the dictionary consisting of all handwritten digits, the digits with same annotation can each be thought of representing a sub-dictionary, and as the proposed algorithm seeks the most economic representation, it is natural that the representation corresponds to picking the representing atoms from the sub-dictionary with greatest affinity with the  digit that represents the data.

In this example, we run the global hybrid IAS algorithm using the parameters $(r^{(1)},\beta^{(1)},\vartheta^{(1)}) = (1,3/2+10^{-4},10^{-5})$, where all components of the vector $\vartheta^{(1)}$ are assumed equal, as sensitivity is not an issue in this example, and  $(r^{(2)},\beta^{(2)}) = (-1,1)$, with the value $\vartheta^{(2)}$ determined from the compatibility condition (\ref{compatibility}). Furthermore, since the digit images are non-negative, after each update step of the pair $(\alpha,\theta)$, we project the image to the positive cone. A theoretical justification of the projection step was given in \cite{CPrSS}. We switch from the first to the second model in the hybrid IAS scheme when either the relative change in $\theta$ with respect to the $\ell_2$-norm falls below $10^{-3}$ or 80 iterations have been completed.  

\begin{figure}
	\centerline{\includegraphics[width=13cm]{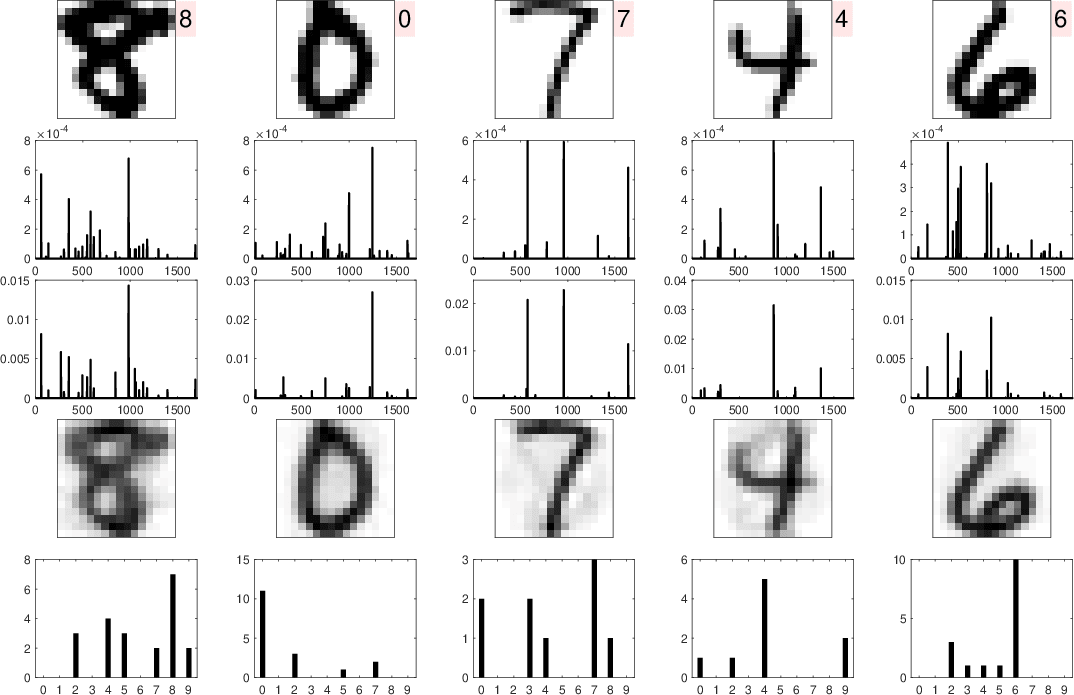}}
	\caption{\label{fig:DL1} Dictionary learning results. The first row shows the test images of the digits to be classified by the dictionary learning algorithm (vector $b$), the true annotation indicated in the figure, the second row the vectors $\theta$ after the IAS iteration with the first hyperprior, and the third row after the iteration with the second hyperprior. The fourth row represents the synthesis $\mW\alpha$ approximating the original digit, and finally, the fifth row gives the histogram of the annotations of the atoms corresponding to coefficients above a threshold $\tau = 0.01$. The annotation is done by majority vote, choosing the largest of the bins. In this example, the standard deviation of the noise representing the mismatch was $\sigma = 0.01$.}
\end{figure}

\begin{figure}
	\centerline{
		\includegraphics[width=13cm]{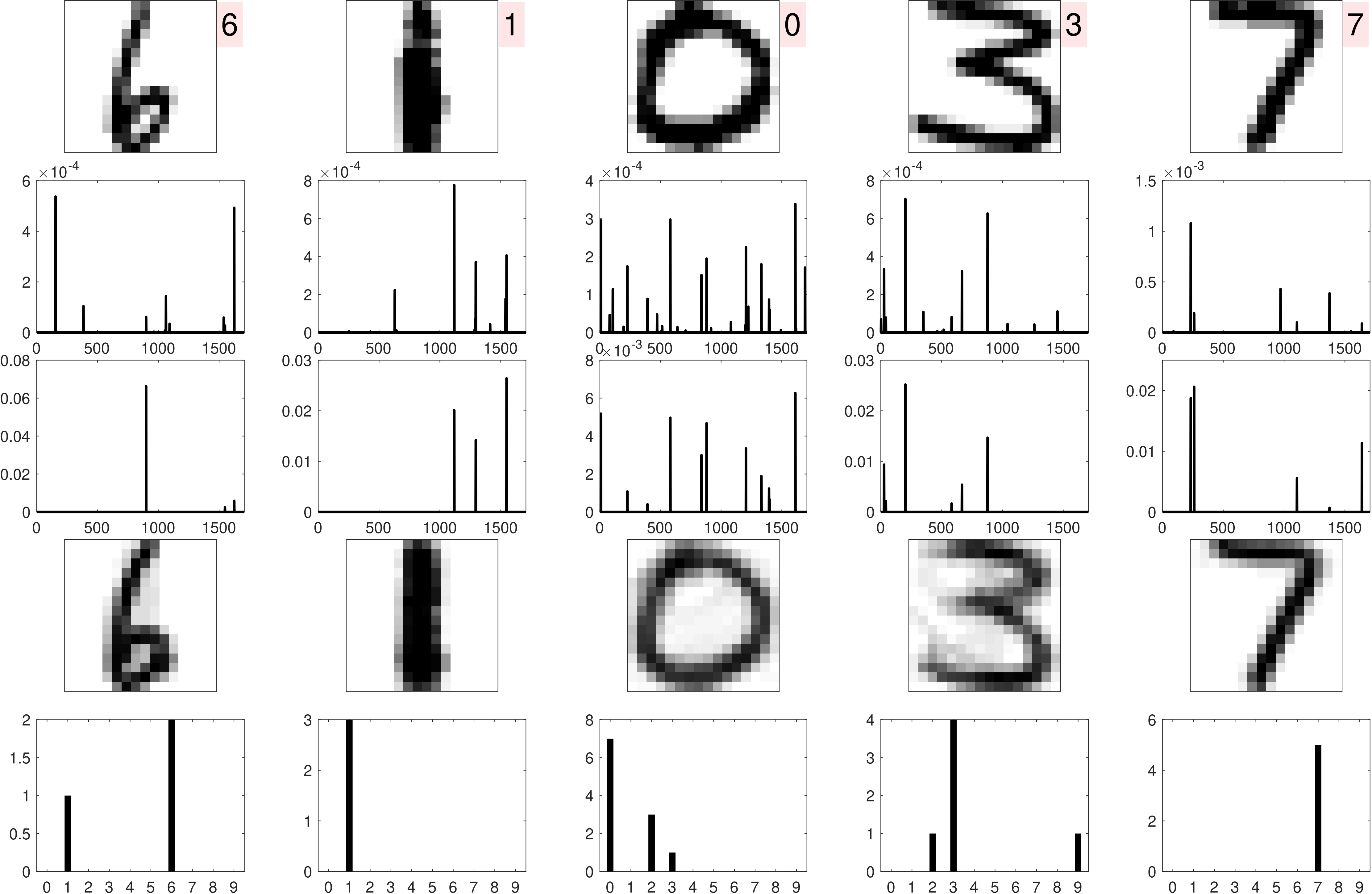}
	}
	\caption{\label{fig:DL2} The rows are as in Figure~\ref{fig:DL1}. In this example, the standard deviation of the noise representing the mismatch was $\sigma = 0.05$. Observe that the approximation becomes sparser.}
\end{figure}

\begin{figure}
	\centerline{
		\includegraphics[width=13cm]{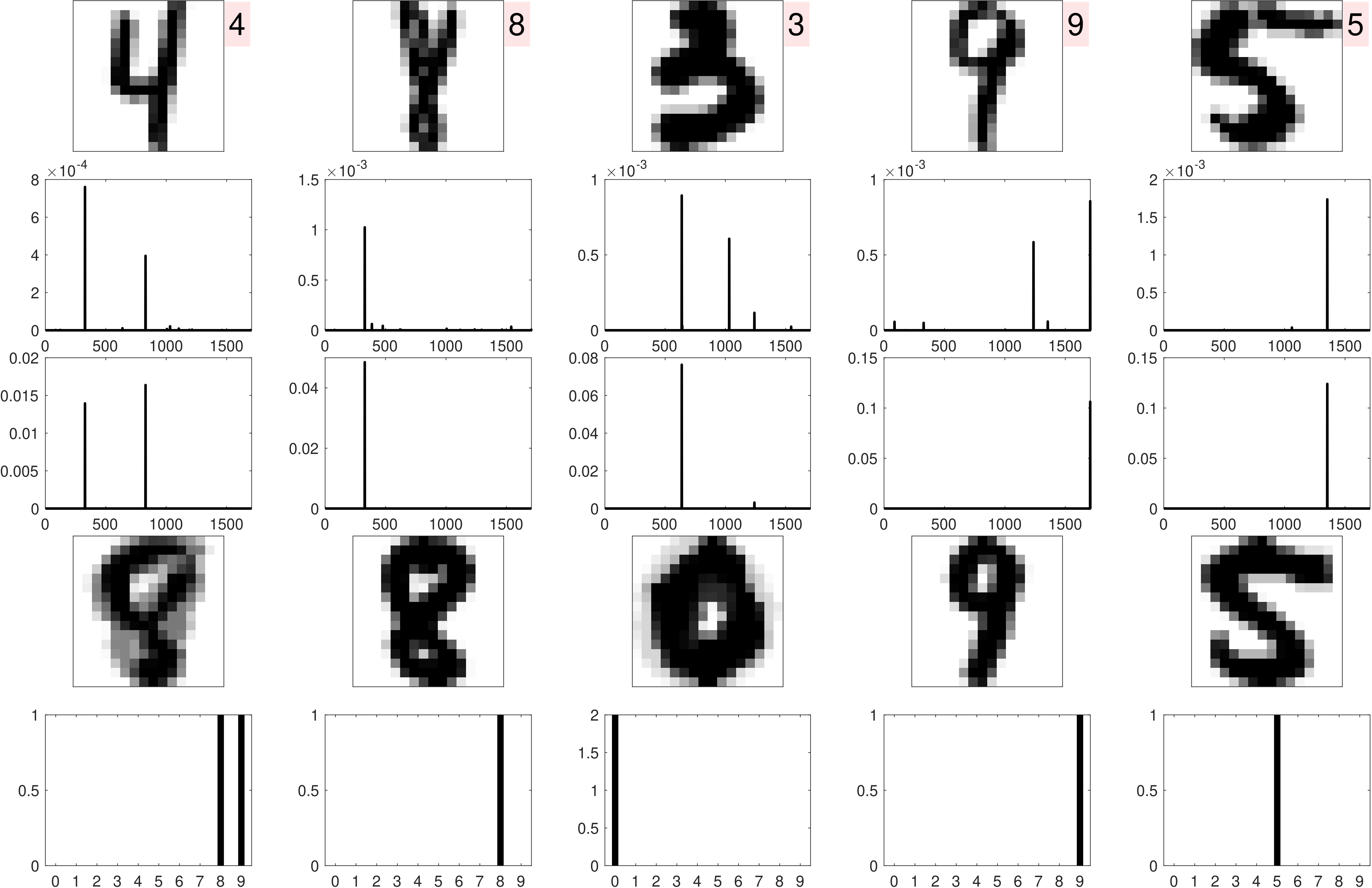}
	}
	\caption{\label{fig:DL3} The rows are as in Figure~\ref{fig:DL1}. In this example, the standard deviation of the noise representing the mismatch was $\sigma = 0.1$. The increased sparsity here is traded with an increased number of misclassifications, such as in the first and the third columns.}
\end{figure}

Figures~\ref{fig:DL1}, \ref{fig:DL2} and \ref{fig:DL3} show the results with different choices of the standard deviation of the likelihood. Observe that here, the noise term $\varepsilon$ represents the discrepancy between the data $b$ and its representation in terms of the dictionary, and can be chosen according to how much fidelity is required. Choosing $\sigma$ large allows a very sparse representation, as the required quality of the approximation is low, however, poor approximation easily leads to a mis-labeling of the digit. On the other hand, decreasing $\sigma$ forces the approximation to be better, and more atoms are required. The labeling can be done using the majority vote principle. In the computed example, Figures~\ref{fig:DL1} and \ref{fig:DL2}, the labeling with majority vote is correct in each case, while in Figure~\ref{fig:DL3} with sparser representation, mislabelings occur.

\section{Conclusions and future work}
\label{sec:concl}
The hierarchical Bayesian framework combined with Krylov subspace iterative solvers for large linear systems is well suited for the design of computationally efficient methods to solve large scale ill posed inverse problems with sparsity constraints. From the point of view of computational efficiency, replacing a whitened Tikhonov-type penalty by a Krylov subspace iteration equipped with early stopping is of crucial importance. The estimate found by this process is not guaranteed to coincide with the MAP estimate, and it was referred to as quasi-MAP estimate in \cite{CPPSV}. The statistical analysis of the early stopping regularization strategy based on Krylov subspace methods is not straightforward, since the estimate depends non-linearly on the data, and it can be seen as an approximate Bayesian computing (ABC) strategy. For further discussion, see, e.g. \cite{SIREV}.
Here we have shown that the framework can be naturally adapted for dealing with overcomplete systems, consisting of,  e.g., combined frames or bases. The approach has significant potential when it may not be known a priori which frame is best suited for representing the unknown, leaving it up to the algorithm to find the most parsimonious representation. In order to avoid that one frame is favored over another, however, it is important that the data are equally sensitive to components in every frame. Fortunately, the sensitivity analysis developed by the authors in \cite{calvetti2019brain,CSS,CPrSS}, provides naturally such scaling. The proposed sensitivity weights are rooted in the very natural Bayesian principle of exchangeability, stating that no set of non-zero components with a given cardinality should be favored over any other. In light of this principle, the scaling guarantees the same explanation power for every sub-frame, so the one leading to most sparse solution is automatically selected. This feature may turn out particularly useful in machine learning, with applications such as MRI fingerprinting (see,e.g., \cite{ma2013magnetic}). 
In \cite{CSS}, a connection between the proposed IAS algorithm and the compressed sensing literature  \cite{candes1} was considered, suggesting that when the forward model guarantees perfect sparse recovery, the IAS algorithm effectively finds a good approximation of it. It is reasonable to believe that the results can be extended to overcomplete dictionaries, for which similar recovery results are known \cite{drip}. 

The methodology developed in this paper has been tested only with Gaussian noise, leading to a quadratic fidelity term in the optimization problem. The IAS framework has been shown to work well with other noise models, e.g., Poisson distributed noise in connection with low dose X-ray tomography and PET, see \cite{BCS}. The applicability of the approach to non Gaussian noise may be very important for its use for dictionary learning problems where the data consist, e.g., of word counts. The extension of the method to large scale problems, different noise models and nonlinear forward models is the next step and will be addressed in separate future contributions.


\begin{thebibliography}{99}
	
	
	\bibitem{BCS}
	\newblock  J.M. Bardsley, D. Calvetti and E. Somersalo  
	\newblock \emph{Hierarchical regularization for edge-preserving reconstruction of PET images},
	\newblock Inverse Problems, \textbf{26}(3), p.035010.
	
	
	\bibitem{boyd2011distributed}
	\newblock  S. Boyd, N. Parikh, E. Chu, B. Peleato and J. Eckstein,  
	\newblock \emph{Distributed optimization and statistical learning via the alternating direction method of multipliers},
	\newblock Foundations and Trends in Machine learning, \textbf{3}(1) (201), 1--122.
	
	
	\bibitem{dict2}
	\newblock  A. M. Bruckstein, D. L. Donoho and M. Elad,  
	\newblock \emph{From Sparse Solutions of Systems of Equations to Sparse Modeling of Signals and Images},
	\newblock SIAM Review, \textbf{51}(1) (2009), 34--81.
	
	
	
	\bibitem{CHPS}
	\newblock  D. Calvetti, H. Hakula, S. Pursiainen and E. Somersalo,  
	\newblock \emph{Conditionally {G}aussian Hypermodels for Cerebral Source Localization},
	\newblock SIAM Journal on Imaging Sciences, \textbf{2}(3) (2009), 879--909.
	
	
	\bibitem{CPPSV}
	\newblock  D. Calvetti, F. Pitolli, J. Prezioso,  E. Somersalo and b. Vantaggi,  
	\newblock \emph{Priorconditioned CGLS-based quasi-MAP estimate, statistical stopping rule, and ranking of priors},
	\newblock SIAM Journal of Scientific Computing, \textbf{39} (2017), S477--S500.	
	
	
	
	\bibitem{calvetti2019brain}
	\newblock  D. Calvetti, A. Pascarella, F. Pitolli, E. Somersalo and B. Vantaggi,  
	\newblock \emph{Brain activity mapping from {MEG} data via a hierarchical {B}ayesian algorithm with automatic depth weighting},
	\newblock Brain topography, \textbf{32}(3) (2019), 363--393.
	
	
	
	\bibitem{SIREV}
	\newblock  D. Calvetti, F. Pitolli,  E. Somersalo and b. Vantaggi,  
	\newblock \emph{Bayes meets Krylov: Statistically inspired preconditioners for CGLS},
	\newblock SIAM Review, \textbf{60} (2018),  429--461.
	
	
	
	
	
	
	\bibitem{CPrS}
	\newblock D. Calvetti, M. Pragliola and E. Somersalo,  
	\newblock \emph{Sparsity promoting hybrid solvers for hierarchical {B}ayesian inverse problems},
	\newblock SIAM Journal on Scientific Computing \textbf{42} (2020), A3761--A3784.
	
	
	
	
	
	
	
	
	\bibitem{CPrSS}
	\newblock D. Calvetti, M. Pragliola, E. Somersalo and A. Strang,  
	\newblock \emph{Sparse reconstructions from few noisy data: analysis of hierarchical {B}ayesian models with generalized gamma hyperpriors},
	\newblock  Inverse Problems, \textbf{36}(2) (2020), p.025010.
	
	
	
	\bibitem{CSS}
	\newblock D. Calvetti, E. Somersalo and A. Strang,  
	\newblock \emph{Hierachical Bayesian models and sparsity:
		$\ell_2$-magic},
	\newblock  Inverse Problems, \textbf{35}(3) (2019), p.035003.
	
	
	
	
	
	
	
	\bibitem{candes1}
	\newblock E. J. {Candes}, J. {Romberg} and T. {Tao},  
	\newblock \emph{Robust uncertainty principles: exact signal reconstruction from highly incomplete frequency information},
	\newblock  IEEE Transactions on Information Theory, \textbf{52}(2) (2006), 489--509.
	
	
	
	
	
	
	\bibitem{drip}
	\newblock E. J. Candes, Y. C. Eldar, D. Needell, D. and P. Randall,  
	\newblock \emph{Compressed sensing with coherent and redundant dictionaries},
	\newblock  Applied and Computational Harmonic Analysis, \textbf{31}(1) (2011), 59--73.
	
	
	
	
	
	\bibitem{chamb}
	\newblock A. Chambolle, M. Holler and T. Pock,  
	\newblock \emph{A Convex Variational Model for Learning Convolutional Image Atoms from Incomplete Data},
	\newblock  Journal of Mathematical Imaging and Vision, \textbf{62} (2020), 417--444.
	
	
	
	\bibitem{dict3}
	\newblock G. Chen and D. Needell, 
	\newblock \emph{Compressed sensing and dictionary learning},
	\newblock  Finite Frame Theory, Proceedings of Symposia in Applied Mathematics, \textbf{73} (2016), 201--241.
	
	
	
	\bibitem{donoho}
	\newblock S. S. Chen, D. L. Donoho, and M. A. Saunders, 
	\newblock \emph{Atomic Decomposition by Basis Pursuit},
	\newblock  SIAM Journal on Scientific Computing, \textbf{20}(1) (1998), 33--61.
	
	
	
	
	
	\bibitem{ma2013magnetic}
	\newblock D. Ma, V. Gulani, N.  Seiberlich, K. Liu, J. L. Sunshine, J. L. Duerk and M. A. Griswold, 
	\newblock \emph{Magnetic resonance fingerprinting},
	\newblock Nature, \textbf{495}(7440) (2013), 187--192.
	
	
	
	\bibitem{mallat}
	\newblock S. G. {Mallat} and Z. Zhang, 
	\newblock \emph{Matching pursuits with time-frequency dictionaries},
	\newblock IEEE Transactions on Signal Processing, \textbf{41}(12) (1993), 3397--3415.
	
	
	
	
	\bibitem{dict1}
	\newblock R. Rubinstein, A. M. {Bruckstein} and M. {Elad}, 
	\newblock \emph{Dictionaries for Sparse Representation Modeling},
	\newblock Proceedings of the IEEE, \textbf{98}(6) (2010), 1045--1057.
	
	
	
	
	
	\bibitem{SFM}
	\newblock J. {Starck}, J. {Fadili} and F. J. {Murtagh}, 
	\newblock \emph{The Undecimated Wavelet Decomposition and its Reconstruction},
	\newblock IEEE Transactions on Image Processing, \textbf{16}(2) (2007), 297--309.
	
	
	
	
	\bibitem{SED}
	\newblock J. L. Starck, M. Elad  and D. Donoho, 
	\newblock \emph{Redundant multiscale transforms and their application for morphological component separation},
	\newblock Advances in Imaging and Electron Physics, \textbf{132} (2004), 287--348.
	
	\bibitem{Vidal}
	\newblock A.F. Vidal, V. De Bortoli, M. Pereyra and A. Durmus, 
	\newblock \emph{Maximum Likelihood Estimation of Regularization Parameters in High-Dimensional Inverse Problems: An Empirical Bayesian Approach Part I: Methodology and Experiments},
	\newblock SIAM Journal on Imaging Sciences, \textbf{13}(4) (2020), 1945--1989.
	
\end{thebibliography}
\end{document}